\documentclass[smallextended]{svjour3}
\pdfoutput=1
\journalname{Journal of Elasticity}
\smartqed  \usepackage{graphicx}
\usepackage{fix-cm}
\usepackage{amsmath}
\usepackage{xcolor}
\usepackage{etoolbox}\newtoggle{COLORWARNINGS}
\newtoggle{NLO}
\newtoggle{DEBUG}
\toggletrue{NLO}
\togglefalse{DEBUG}
\toggletrue{COLORWARNINGS}

\iftoggle{DEBUG}
{
\overfullrule=1mm }{
}
\iftoggle{COLORWARNINGS}
{

}{

}

\usepackage{amssymb}\usepackage{siunitx}\usepackage{suffix}\usepackage[hyphens]{url}\usepackage{booktabs} \togglefalse{NLO}

\newcommand{\PD}[2]{\frac{\partial #1}{\partial #2}}

\newcommand{\PDDX}[3]{\frac{\partial^2 #1}{\partial #2 \partial #3}}
\newcommand{\D}[2]{\frac{d #1}{d #2}}
\newcommand{\ABS}[1]{\left|#1\right|}
\newcommand{\SYMPART}[1]{\frac{1}{2}\left( #1 + (#1)^T \right)}
\newcommand{\NORM}[1]{\|{#1}\|}

\makeatletter
\newcommand{\TRANSPOSE}{\@ifstar
                     \TRANSPOSESTAR                     \TRANSPOSENOSTAR}

\newcommand{\TRANSPOSENOSTAR}[1]{\left({#1}\right)^T}
\newcommand{\TRANSPOSESTAR}[1]{{#1}^T}
\makeatother

\DeclareMathOperator{\graaad}{\nabla}
\newcommand{\grad}[1]{\graaad\!#1}

\DeclareMathOperator{\divergence}{div}
\newcommand{\diver}[1]{\divergence #1}

\DeclareMathOperator{\tr}{tr}

\newcommand{\tenzor}[1]{\boldsymbol{#1}}

\newcommand{\Esymbol}{\varepsilon}

\newcommand{\Tsymbol}{T}

\newcommand{\Ex}[1]{\Esymbol_{#1}}

\newcommand{\Tx}[1]{\Tsymbol_{#1}}

\newcommand{\Ux}[1]{u_{#1}}

\newcommand{\at}[2]{\left.#1\right|_{#2}}

\newcommand{\E}{\tenzor{\Esymbol}}
\newcommand{\T}{\tenzor{\Tsymbol}}
\newcommand{\XU}{\vec{u}}
\newcommand{\IET}{\tenzor{G}}

\newcommand{\INTEGRAL}[3]{\int#1 \! #2 \, \mathrm{d} #3}

\newcommand{\RN}[1]{  \textup{\uppercase\expandafter{\romannumeral#1}}}

\newcommand{\defref}[1]{Definition~\ref{#1}}
\WithSuffix\newcommand\defref*[1]{Definition~\ref{#1}, (p.~\pageref{#1})}

\newcommand{\theoref}[1]{Theorem~\ref{#1}}
\WithSuffix\newcommand\theoref*[1]{Theorem~\ref{#1}, (p.~\pageref{#1})}

\newcommand{\lemref}[1]{Lemma~\ref{#1}}
\WithSuffix\newcommand\lemref*[1]{Lemma~\ref{#1}, (p.~\pageref{#1})}

\newcommand{\secref}[1]{Section~\ref{#1}}
\WithSuffix\newcommand\secref*[1]{Section~\ref{#1}, (p.~\pageref{#1})}

\newcommand{\obsref}[1]{Observation~\ref{#1}}
\WithSuffix\newcommand\obsref*[1]{Observation~\ref{#1}, (p.~\pageref{#1})}

\newcommand{\figref}[1]{Figure~\ref{#1}}
\WithSuffix\newcommand\figref*[1]{Figure~\ref{#1}, (p.~\pageref{#1})}

\newcommand{\tabref}[1]{Table~\ref{#1}}
\WithSuffix\newcommand\tabref*[1]{Table~\ref{#1}, (p.~\pageref{#1})}

\begin{document}
\title{The state of stress and strain adjacent to notches in a new class of nonlinear elastic bodies.\thanks{Josef M\'alek thanks the Czech Science Foundation for support through the project 18-12719S.}
\thanks{K.R.~Rajagopal thanks the Office of Naval Research for support of this work.}
}
\titlerunning{Computer simulations of power law solids}

\author{Vojt\v ech Kulvait \and
Josef M\' alek \and 
K.R.~Rajagopal
\thanks{Dedicated to the memory of Professor Walter Noll.}
}

\institute{Vojt\v ech Kulvait \and Josef M\' alek \at Charles University, Faculty of Mathematics and Physics, Mathematical Institute, Sokolovsk\' a 83, 186 75, Prague, Czech Republic 
\and 
K.R. Rajagopal \at Texas A\&M University, Department of Mechanical Engineering, College Station, TX, 77843, USA}

\date{Received: TBD / Accepted: TBD}

\maketitle

\begin{abstract}
In this paper we study the deformation of a body with a notch subject to an anti-plane state of stress within the context of a new class of elastic models. These models stem as approximations of constitutive response functions for an elastic body that is defined within the context of an implicit constitutive relation between the stress and the deformation gradient. Gum metal and many metallic alloys are described well by such constitutive relations. We consider the state of anti-plane stress of a body with a smoothened V-notch within the context of constitutive relations for the linearized strain in terms of a power-law for the stretch. The problem is solved numerically and the convergence and the stability of the solution is studied.

\end{abstract}
 \section{Introduction}
Walter Noll \cite{KRR_NOLL1955,KRR_NOLL1957,KRR_NOLL1958} introduced the concept of a Simple Material\footnote{Noll \cite{KRR_NOLL1972} later generalized the concept of his definition of a Simple Material. We shall not get into a discussion of the same here.} whose constitutive representation subsumes many classical constitutive expressions for elastic and viscoelastic solids and viscous and viscoelastic fluids, thus providing a framework within which one could study the response characteristics of a large class of materials. This notwithstanding, since the Cauchy stress in a Simple Material depends on the history of the density and the deformation gradient it precludes the possibility that the histories of the stress and the deformation gradient could be related by an implicit constitutive relation. Recently, Rajagopal \cite{KRR_2003,KRR_2006} has provided a compelling rationale for considering such implicit constitutive relations to describe the response of both fluids and solids. In general, one can have an implicit relationship between the history of the density, stress, deformation gradient, and possibly other relevant physical variables such as the temperature, electrical and magnetic fields, etc\footnote{Rajagopal \cite{C_Rajagopal2015} classsified the material symmetry possessed by the sub-classes of bodies whose histories of the stress, deformation gradient, density, etc., are given in terms of an implicit constitituve relations.}.

As algebraic implicit constitutive relations are somewhat recent, we provide some related references to studies concerning the response of bodies described by implicit constitutive relations. Pr{\r{u}}{\v{s}}a and Rajagopal \cite{KRR_PRUSA2012} generalizing the approach of Coleman and Noll \cite{KRR_COLEMAN1960} were able to obtain constitutive response relations of the type due to Maxwell, Oldroyd, Rivlin and Ericksen, that are approximations that hold within the context of retarded motions. Perl{\'{a}}cov{\'{a}} and Pr{\r{u}}{\v{s}}a \cite{D_Perlacova2015}, building on the earlier work of LeRoux and Rajagopal \cite{E_Roux2013}, generalizing it and using algebraic implicit constitutive relations (that is constitutive relations wherein only the stress and the symmetric part of the velocity gradient, and not any of their time derivatives, appear) used them describe the response of colloids and suspensions. M\'{a}lek et al. \cite{F_Malek2010} developed a generalization of the Navier-Stokes constitutive relation and introduced stress power-law fluids and Rajagopal \cite{G_Rajagopal2013} showed that had one started by expressing the velocity gradient as a function of the stress than vice-versa one would not arrive at erroneous assumptions such as the Stokes assumption in fluid mechanics. 

Rigorous mathematical analyses for a subclass of fluids described by implicit constitutive relations (including the activated fluids such as the Bingham fluid or activated Euler fluids) have been carried out for steady flows by M\'{a}lek et al.  \cite{Malek_ETNA,BGMS_ACV} and for unsteady flows by Bul\'\i{}\v{c}ek et al. \cite{H_Bulicek2012,I_Buliccek}; see also a recent study by Blechta et al. \cite{BMR_2018} where the systematic classification of fluids described by implicit constitutive relation are carried out and the mathematical analysis of steady/unsteady flows of activated Euler fluids subject to various types of boundary conditions is presented. Earlier, the analysis of unsteady flows of Bingham fluids subject to boundary conditions described by implicit relations is considered in \cite{BM_2016} and flows of heat conducting fluids are analysed in \cite{BulMal_Vietnam,MZ_2018}. The convergence of numerical schemes proposed to approximate flows of incompressible fluids towards the (weak) solution of the original PDE problems is established in \cite{DKS_2013,KS_2016,ST_2018}. PDE analysis for a class of solids with bounded linearized strain is developed in \cite{BMRW2015,BMS_2015,Beck_2017}, see also the survey paper \cite{bulicek_elastic_2014}, while a study regarding the properties of finite element approximations is presented in \cite{BGS_2018} (see also \cite{kulvait_anti_2012,Montero2016} for relevant results of computer simulations). Large data and long time PDE analysis of unsteady motions of generalized Kelvin-Voigt solids is developed in \cite{BMR_2012}; see also the study \cite{BGMS_M3AS}. Finally, the analysis of the equations governing the motion of a special class of compressible fluids with bounded divergence is developed in \cite{Feireisl_Liao_Malek2015}, see also \cite{M_KRR_Comp} for derivation of these models using a thermodynamical approach.

The response of solids described by implicit constitutive relations has been studied in the case of polymers by Rajagopal and Saccomandi \cite{J_Rajagopal2009}, inelastic solids by Rajagopal and Srinivasa \cite{K_Rajagopal2015} and compressible fluids in M\'{a}lek and Rajagopal \cite{M_KRR_Comp}. The response of electroelastic bodies described by implicit constitutive relations have been studied in the articles by Bustamante and Rajagopal \cite{L_Bustamante2012}, \cite{M_Bustamante2013}, and the response of magnetoelastic bodies described by implicit constitutive relations has been analysed by Bustamante and Rajagopal \cite{N_Bustamante2015}.  A review of implicit constitutive theories can be found in the articles by Rajagopal and Saccomandi \cite{O_Rajagopal2016} and M\'{a}lek and Pr{\r{u}}{\v{s}}a \cite{MP_2018}.

When one confines attention to elastic solids, a special sub-class of the implicit constitutive relations that is relevant is the class of bodies whose response is described by an implicit relationship between the stress and the deformation gradient. Rajagopal and Srinivasa \cite{A_Rajagopal2006}, \cite{B_Rajagopal2009} have provided a thermodynamic basis for implicit theories of elasticity. Cauchy elastic bodies wherein the stress is expressed as an explicit function of the deformation gradient are a special subclass of the above class of implicit relations. Another special subclass is defined by an explicit expression for the Cauchy-Green tensor $\tenzor{B}$ in terms of the Cauchy stress $\T$. Truesdell and Moon \cite{KRR_TRUESDELL1975} studied conditions that guaranteed semi-invertibility of isotropic functions and thus is relevant to the conditions under which the relation between the Cauchy-Green tensor $\tenzor{B}$ and the Cauchy stress $\T$ is invertible. But not all functions of the Cauchy-Green tensor $\tenzor{B}$ as a function of the Cauchy stress $\T$ are invertible. The class of response functions that express the Cauchy-Green tensor $\tenzor{B}$ as a function of the Cauchy stress $\T$ includes response functions that are not Cauchy elastic (Cauchy \cite{KRR_CAUCHY1822}). Moreover, Truesdell and Moon \cite{KRR_TRUESDELL1975} were not interested in issues of causality that would lead one to conclude that implicit models that relate the stress to kinematical quantities have a sounder philosophical underpinning. Moreover, since the study of Truesdell and Moon \cite{KRR_TRUESDELL1975} was confined to isotropic response, it does not address any of the models that describe the relationship between the Cauchy-Green tensor $\tenzor{C}$ and the Cauchy stress in the case of anisotropic response of elastic bodies.  

There are several advantages to the use of implicit models in describing elastic response, one of them being the fact that the use of the linearization based on the displacement gradient being small leading to a nonlinear relationship between the linearized strain and the stress (see Rajagopal \cite{KRR_2007,rajagopal_conspectus_2011}). Such responses are exhibited by a variety of metallic alloys as well as materials like concrete and Devindran et al. \cite{devendiran_thermodynamically_2016} used such implicit models to describe the response of several metallic alloys. Also, a sub-class of these models, those that have limited linearized strain\footnote{Bul\'{i}\v{c}ek et al. \cite{bulicek_elastic_2014} study several mathematical aspects concerning strain limiting bodies.} are particularly well suited for the study of problems such as fracture and the state of stresses and strains at notches within the context of brittle materials (see Rajagopal and Walton \cite{rajagopal_modeling_2011}, Gou et al. \cite{gou_modeling_2015}, Kulvait et al. \cite{kulvait_anti_2012}, Ortiz et al. \cite{ortiz_numerical_2012,ortiz_numerical_2014}, Montero et al. \cite{Montero2016}, Itou et al. \cite{itou_nonlinear_2016}) and unlike the classical linearized theory of elasticity the strains are not singular. Of course, this is not surprising as the constitutive relation has built into it the concept of bounded strains.
The problem of anti-plane state of stress in a body with a V-notch that is pointed or sharply radiused has been studied by Zappolorto et al. \cite{XZappalorto2016} in which they obtain closed form solutions in the case of strain limiting models, and Rajagopal and Zappalorto \cite{Q_Rajagopal2018} have studied the state of stress and  strain adjacent to a crack tip in the case of a non-monotonic relationship between the stress and the strain.

Instead of strain limiting theories, one could have constitutive relations wherein the relationship between the linearized strain and the stress is given by a power-law. Such models have been found to be very useful in describing the response of Gum metal and several Titanium alloys (see Rajagopal \cite{rajagopal_nonlinear_2014}, Kulvait et al. \cite{kulvait_modeling_2017}). In such models, unlike strain limiting models, one finds that appropriate values of the power-law index, the strains grows as the stress grows, but much slower than that for a linearized elastic model. It would be interesting to examine the state of strain and stress at notches which have been smoothed (see \figref{fg5:computationalDomainVC}) and it is towards the investigation of this question that this study is aimed. Two parameters, an angle $\alpha$ and a radius $r_c$ determine the nature of the smoothening that is carried out. We study the anti-plane state of stress for a rectangle with a smoothened notch within the context of different materials. The material characterization for a variety of Titanium alloys has been carried out by Kulvait et al. \cite{kulvait_modeling_2017} and we use these material parameters to characterize the different alloys that comprise the rectangular body with a notch that is subjected to anti-plane stress in the study.

In the present study we are interested in the class of constitutive relations, where the spherical (usually referred to as "hydrostatic part" or "isotropic part", but neither of these terminologies seem appropriate) and deviatoric parts of the deformation response are modeled through
\begin{equation}
\begin{aligned}
\tr{\E} &= \sigma_1(\tr{\T}) \tr{\T},\\
\E^d &= \sigma_2(|\T^d|) \T^d,\label{eq2:deviatoric}
\end{aligned}
\end{equation}
where $\sigma_1$ and $\sigma_2$ are scalar functions, with $\sigma_1(0)=0$ and $\sigma_2(0)=0$, $\T$ is a Cauchy stress, $\E$ is the linearized strain tensor, $\tr \E$ and $\E^d$ denote the trace and deviatoric part of the tensor $\E$. Of course, one could provide the form for the linearized strain in terms of the Cauchy stress that implies the assumed forms for spherical and deviatoric parts, but we have chosen to express it as above so that the structure of the two parts of the linearized strain are transparent. In particular we are interested in power-law models, wherein
\begin{equation}
\begin{aligned}
\sigma_1(\tr{\T}) &= \frac{1}{3 K}\left(\frac{\tau_K^2 + |\tr{\T}|^2}{\tau_K^2}\right)^{\frac{s'-2}{2}},\\
\sigma_2(|\T^d|) &= \frac{1}{2 \mu} \left(\frac{\tau_\mu^2 + |\T^d|^2}{\tau_\mu^2}\right)^{\frac{q'-2}{2}}, \label{eq2:powerlaw1}
\end{aligned}
\end{equation}
parameters $q' \in (1, \infty)$, $s' \in (1, \infty)$, $K>0$, $\tau_K>0$, $\mu>0$,  $\tau_\mu>0$ are material moduli. One should not think of expressions obtained by inverting \eqref{eq2:deviatoric}, with \eqref{eq2:powerlaw1} inserted in the appropriate place, as constitutive models. This is because these models have the meaning of approximations only within the context of the expressions for the Cauchy-Green tensor $\tenzor{B}$ in terms of the stress $\T$. A detailed discussion concerning the relevant issues can be found in Rajagopal \cite{KRR_2017}.

Recently, it has been shown that the model \eqref{eq2:powerlaw1} is capable of capturing the nonlinear response characteristics of titanium alloys in the small strain range, see \cite{kulvait_modeling_2017}.  We have identified the parameters of the power-law model \eqref{eq2:powerlaw1}, for four different beta phase titanium alloys. Based on these identified parameters, we have studied the behavior of these alloys when subject to anti-plane shear stress. 

\begin{figure}[h]
\centering \includegraphics{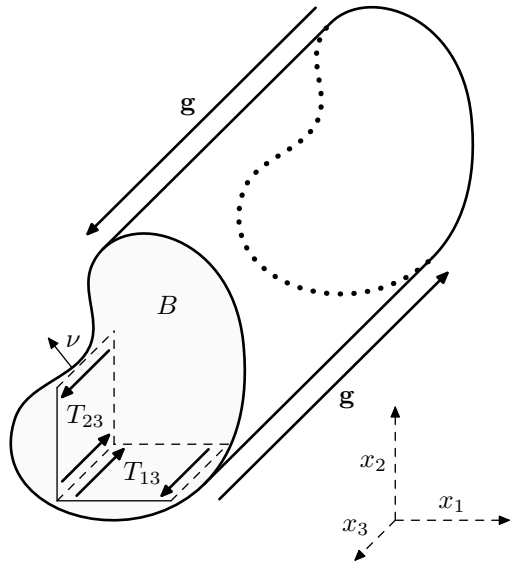}
    \caption[Anti-plane stress setting]{Anti-plane stress setting. The only nonzero components of the stress tensor are $T_{13}$ and $T_{23}$.}
    \label{img5:antiPlaneMotivation}
\end{figure}

This study is devoted to the anti-plane stress of a body containing a smoothed V-notch, which we refer to as the VC notch. As it is an anti-plane stress problem, we shall assume that the stress and the displacement field depend only on plane coordinates ($x_1$ and $x_2$) and the only nonzero components of the stress tensor are $\Tx{13}$ and $\Tx{23}$. The boundary tractions on a body that is subject to the anti-plane stress setting is in a direction perpendicular to the ($x_1$, $x_2$) plane, see Figure \ref{img5:antiPlaneMotivation}. These assumptions allow us to rephrase the problem of finding unknown stress and displacement of the body as a variational problem on a discrete space of finite elements to find an unknown Airy's stress function. 

In \cite{kulvait_anti_2012} we have presented numerical study of the behavior of the strain-limiting model in the anti-plane stress setting. Here, we are interested in the behavior of the class of materials defined through \eqref{eq2:deviatoric} and \eqref{eq2:powerlaw1} in the geometry of a square plate with a smoothened V-shaped notch.

The software implementation for solving these problems in Python can be found in supplementary materials. Solver of the discretized problem on the finite element space uses the damped Newton method and utilizes FEniCS software library, see \cite{logg_automated_2012}.
 \section{Boundary value problem}
The fundamental boundary value problem of elastostatics can be formulated as follows.
\begin{definition}[Problem (P)]
Let $\Omega \subset \mathbb{R}^n$, $n = 2$ or $3$ be an open, bounded, connected set with the boundary $\partial \Omega$ consisting of two smooth disjoint parts $\Gamma_D$ and $\Gamma_N$ such that $\partial \Omega = \overline{\Gamma_D \cup \Gamma_N}$, where $\vec{\nu} = \vec{\nu}(\vec{x})$ denotes the unit outward normal vector at $\vec{x} \in \partial \Omega$. Let $\IET : \mathbb{R}^{n \times n}_{sym} \to \mathbb{R}^{n \times n}_{sym}$, $\vec{u_0} : \Omega \to \mathbb{R}^n$, $\vec{f} : \Omega \to \mathbb{R}^n$ and $\vec{g}: \Gamma_N \to \mathbb{R}^n$ be given. We say that the pair of functions $(\XU, \T) : \Omega \to \mathbb{R}^n \times \mathbb{R}^{n \times n}_{sym}$ solves the Problem~(P) if
\begin{subequations}
\begin{align}
-\diver{\T} &= \vec{f}&&\mbox{ in } \Omega, \label{eq1:equilibrium}\\
\E &= \IET(\T)&&\mbox{ in } \Omega, \label{eq3:constrel}\\
\XU &= \vec{u_0} &&\mbox{ on } \Gamma_D,\label{eq1:dir}\\
\T \vec{\nu} &= \vec{g}&& \mbox{ on } \Gamma_N,\label{eq1:neu}
\end{align}
\label{eq1:analysisbvp}
\end{subequations}
where
\begin{equation}
\E = \SYMPART{\grad{\XU}}. \nonumber
\end{equation}
\label{def01:bvp}
\end{definition}
The above problem describes the state of the elastic body occupying the set $\Omega$ at equilibrium characterized by the equation \eqref{eq1:equilibrium}, where $\T$ is the stress tensor and $\vec{f}$ stands for external body forces. The body obeys the nonlinear constitutive equation \eqref{eq3:constrel} that relates the linearized strain $\E$, which is a symmetric gradient of the displacement $\XU$, and the stress tensor $\T$ by means of a nonlinear function $\IET$. The equation \eqref{eq1:dir} prescribes the displacement $\vec{u_0}$ on the Dirichlet part of the boundary $\Gamma_D$, and the equation \eqref{eq1:neu} prescribes the boundary traction $\vec{g}$ on the Neumann part of the boundary $\Gamma_N$. The state of the body is subject to an additional assumption that the square of the norm of the displacement gradient can be neglected in comparison to the norm of the displacement gradient itself.

The geometry and the type of the deformation of anti-plane stress, see Figure \ref{img5:antiPlaneMotivation}, allows us to assume that $\Omega = B \times \mathbb{R}$, $\XU = (0,0,u(x_1,x_2))$, $\partial B = \Gamma_N$, $\vec{g} = (0,0,g)$ and $\vec{f} = \tenzor{0}$. The stress tensor $\T$ has two nontrivial components and can be represented by vector $\T_v = (\Tx{13}, \Tx{23})$. The trace and the norm of the deviatoric part of the stress tensor can be expressed as
\begin{equation}
\tr{(\T)} = 0, \quad |\T^d| = \sqrt{2} |\T_v|.\nonumber
\end{equation}
We simplify the model \eqref{eq2:powerlaw1} under the assumption that we have a anti-plane stress problem on hand. Then the strain tensor $\E$ has also only two nontrivial components $\Ex{13}$ and $\Ex{23}$ and can be represented by the vector $\E_v = (\Ex{13}, \Ex{23})$. Constitutive relation \eqref{eq2:deviatoric} reduces to the form
\begin{equation}
\E_v = \sigma_2(\sqrt{2} |\T_v|) \T_v \equiv \sigma( |\T_v|) \T_v,\nonumber
\end{equation}
and Problem (P) can be reformulated as:
\begin{definition}[Problem (P) in anti-plane stress setting]
\label{d05:antiPlaneStress}
Let $\Omega = B \times \mathbb{R}$, where $B \subset \mathbb{R}^2$ is an open, bounded, simply connected set. Let $\sigma : \mathbb{R} \to \mathbb{R}$ represent the constitutive response of the material and let $g : \partial B \to \mathbb{R}$ be a given function.  We say that a pair of functions $(\E_v, \T_v)$ is the solution of Problem (P) in the anti-plane stress setting \footnote{Since we consider constitutive relations of the type \eqref{eq4:constrel}, the anti-plane stress state is equivalent to the classical definition of anti-plane strain.} when the following is true
\begin{subequations}
\begin{align}
-\PD{\Tx{13}}{x_1} - \PD{\Tx{23}}{x_2} &= 0 &&\mbox{in } B,\label{eq4:equilibrium}\\
\Ex{13} = \sigma(\ABS{\T_v}) \Tx{13}, \quad \Ex{23} &= \sigma(\ABS{\T_v}) \Tx{23}&&\mbox{in } B,\label{eq4:constrel}\\
\Tx{13}\nu_1 + \Tx{23}\nu_2 &= g &&\mbox{on }\partial B\label{eq4:boundarycondition}.
\end{align}
\label{eq5:bvp}
\end{subequations}
\end{definition}
\vspace{-0.5cm}
In the remainder of this chapter, we use the following particular forms of the constitutive function $\sigma$. For the power-law model \eqref{eq2:powerlaw1}, we have the response
\begin{equation}
\sigma(\ABS{\T_v})  = \frac{1}{2 \mu} \left( \frac{\tau_0^2 + 3\ABS{\T_v}^2}{\tau_0^2} \right)^{\frac{q'-2}{2}} .\label{eq5:powerlaw}
\end{equation}
The linear Hooke's law is characterized by the response \eqref{eq5:powerlaw}, when $q' = 2$, that is
\begin{equation}
\sigma(\ABS{\T_v}) = \frac{1}{2\mu_L}. \label{eq5:hookelaw}
\end{equation}
In the case of strain-limiting model, see \cite{kulvait_anti_2012} we have
\begin{equation}
\sigma(\ABS{\T_v}) = \frac{\tau_\mu}{2 \mu_l \left( \tau_\mu^{\,a} + (\sqrt{2}\ABS{\T_v})^{\,a}\right)^{1/a}}. \label{eq:strainlimitingmodel}
\end{equation}

\subsection{Compatibility conditions}
The Saint-Venant compatibility conditions for $\E$, are reduced to the two nontrivial equations
\begin{equation}
\PDDX{\Ex{13}}{x_1}{x_2} - \PDDX{\Ex{23}}{x_1}{x_1} = 0, \quad \PDDX{\Ex{23}}{x_1}{x_2} - \PDDX{\Ex{13}}{x_2}{x_2} = 0, \label{eq05:compatantiplane}
\end{equation}
which implies the existence of a constant $C$ such that
\begin{equation}
\PD{\Ex{13}}{x_2} - \PD{\Ex{23}}{x_1} = C. \label{eq5:compatantiplane2}
\end{equation}
Expressing \eqref{eq5:compatantiplane2} in terms of $\XU$, we obtain
\begin{equation}
\PDDX{\Ux{1}}{x_3}{x_2} - \PDDX{\Ux{3}}{x_1}{x_2}  - \PDDX{\Ux{2}}{x_3}{x_1} + \PDDX{\Ux{3}}{x_2}{x_1} = C. \label{eq5:compatantiplane3}
\end{equation}
Since  $\XU = (0,0,u(x_1,x_2))$, we conclude that $C=0$ and thus \eqref{eq05:compatantiplane} leads to
\begin{equation}
\PD{\Ex{13}}{x_2} - \PD{\Ex{23}}{x_1} = 0. \label{eq5:compatantiplane4}
\end{equation}
The condition \eqref{eq5:compatantiplane4} is the necessary and sufficient condition for the existence of the displacement $u(x_1,x_2)$ fulfilling
\begin{equation}
\Ex{13} = \frac{1}{2} \PD{u}{x_1}, \quad \Ex{23} = \frac{1}{2} \PD{u}{x_2}. \nonumber
\end{equation}

\subsection{Airy's function}

Airy's stress function $A$ is a scalar function, whose derivatives are components of $\T_v$ defined through
\begin{equation}
\Tx{v}^1 = \Tx{13} = \PD{A}{x_2}, \quad \Tx{v}^2 = \Tx{23} = -\PD{A}{x_1}. \label{eq5:airysfuncion}
\end{equation}
Every stress field of the type \eqref{eq5:airysfuncion} fulfills the equilibrium equation \eqref{eq4:equilibrium}. From \eqref{eq5:airysfuncion}, we have that $|\T_v| = |\grad A|$, and thus $\sigma(|\T_v|) = \sigma(|\grad A|)$. Substituting the constitutive relation \eqref{eq4:constrel} into the compatibility condition \eqref{eq5:compatantiplane4} leads to
\begin{equation}
-\PD{}{x_1} \left( \sigma(\ABS{\nabla A})\PD{A}{x_1} \right) - \PD{}{x_2} \left( \sigma(\ABS{\nabla A})\PD{A}{x_2} \right) = 0 \quad \mbox{in }B. \label{eq4:antiplaneeq}
\end{equation}
The boundary condition \eqref{eq4:boundarycondition} takes the form 
\begin{equation}
\PD{A}{x_2} \nu_1 - \PD{A}{x_1} \nu_2 = \PD{A}{x_1} t_1 +  \PD{A}{x_2} t_2 =  g \qquad \mbox{on } \partial B, \label{eq4:antiplanebcinit}
\end{equation}
where $\vec{t} = (-\nu_2, \nu_1)$ represents the tangential vector to the boundary $\partial B$.
We parametrize the boundary by the counterclockwise oriented closed curve $\vec{\xi} = (\xi_1, \xi_2)$, $\vec{\xi}(0) = \vec{\xi}(h)$, such that $\vec{\xi}([0,h]) = \partial B$ and
\begin{equation}
t_1 (\vec{\xi}(b)) = \xi_1'(b), \quad t_2 (\vec{\xi}(b)) = \xi_2'(b). \label{eq5:vector}
\end{equation}
Substituting \eqref{eq5:vector} into \eqref{eq4:antiplanebcinit} and using the chain rule, we conclude that $g$ is equal to the tangential derivative of $A$. We have that
\begin{equation}
g(\vec{\xi}(b)) = \nabla\!A \cdot \vec{\xi}'(b) = \D{A(\vec{\xi}(b))}{b}, \quad b \in [0,h). \label{eq4:reintegrationbc}
\end{equation}
Integrating \eqref{eq4:reintegrationbc} along the boundary, we obtain the Dirichlet boundary condition for A
\begin{equation}
A(x) = A_0^g(\vec{\xi}(b)) = A(\vec{\xi}(0)) + \INTEGRAL{_{0}^b}{g(\vec{\xi}(s))}{s}, \quad b \in [0,h), x =  \vec{\xi}(b) \in \partial B. \label{eq4:findingbc}
\end{equation}
By introducing Airy's stress function via \eqref{eq5:airysfuncion}, the boundary value problem \eqref{eq5:bvp} takes the form summarized in the following definition.

\begin{definition}[Airy's function for solving Problem (P) in anti-plane stress setting]
\label{d05:antiPlaneStressAiryFunction}
Let $B \subset \mathbb{R}^2$  be an open, bounded, simply connected set and the functions $(g, \sigma)$ fulfill the assumptions of Definition \ref{d05:antiPlaneStress}. We say that the function $A$ is the Airy's function solving Problem (P) in the anti-plane stress setting if the following is true
\begin{equation}
\begin{aligned}
-\diver \left( \sigma(\ABS{\grad A})\grad{A} \right) &= 0 \quad &&\mbox{in }B,\\
A(\vec{x}) &= A_0^g(\vec{x}) \qquad &&\mbox{on } \partial B,
\end{aligned}
\label{eq5:airybvp}
\end{equation}
where $A_0^g$  is the function introduced in \eqref{eq4:findingbc}.
\end{definition}
\vspace{-0.5cm}

 \section{Finite element simulations}
\label{sec05:setting}
In this section, we delineate the setting that we follow when performing computer simulations. First, we depict the geometries of the computational domains and describe the boundary conditions that we use. Then, we define the variational formulation of the problem on finite dimensional spaces. Finally, we tabulate the parameters of the models of the titanium alloys that we simulate.

\subsection{Computational domain}
\label{sec05:cd}
\begin{figure}[h]
\begin{center}
\includegraphics[height=6.7cm]{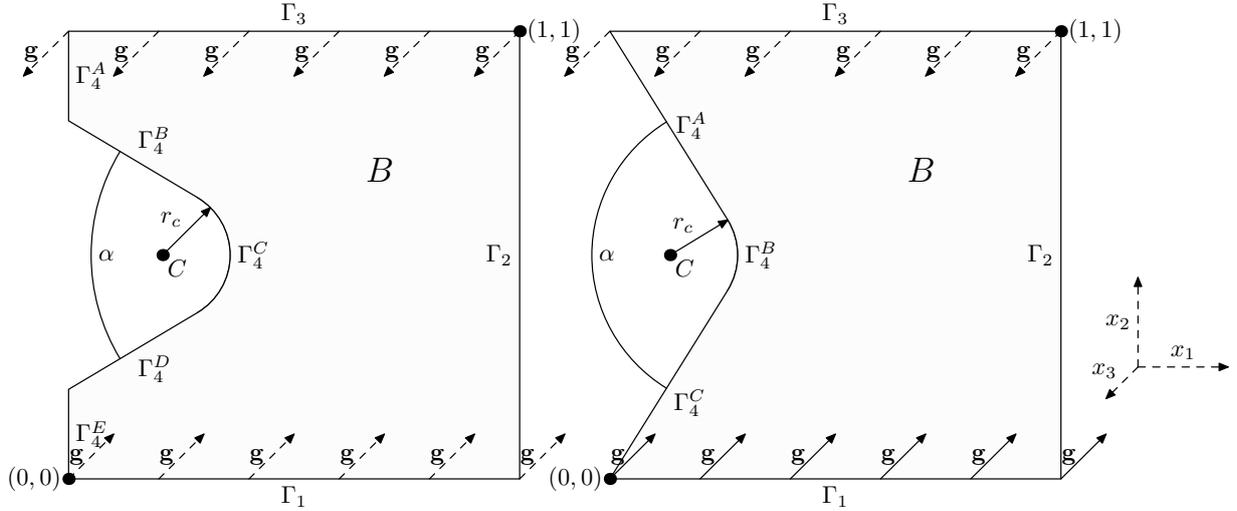}
\end{center}
\caption[VC-geometry of a computational domain]{VC-geometry of a computational domain. On the left, the domain is depicted when $\alpha < \pi/2$. On the right, when $ \pi/2 \leq \alpha < \pi$.}
\label{fg5:computationalDomainVC}
\end{figure}

We use the geometry of the square with a smoothened V-notch, which we shall refer to as the the VC-geometry, see Figure \ref{fg5:computationalDomainVC}. The tip of the V notch is smoothened by the circle arc of radius $r_c$ that is tangent to the V-notch.  The geometry is parametrized by angle $\alpha$ and the radius of arc $r_c$. When $\alpha < \pi / 2$, then the left part of the boundary consists of five parts $\Gamma_4 = \Gamma_4^A \cup \Gamma_4^B \cup \Gamma_4^C \cup \Gamma_4^D \cup \Gamma_4^E$. There is a restriction that the arc fits inside the unit square, which imposes the condition
\begin{equation}
0 < r_c < \frac{0.5 \sin(\frac{\alpha}{2})}{\cos^2(\frac{\alpha}{2})}.\nonumber
\end{equation}
The center of the arc is at the point 
\begin{equation}
C = \left( 0.5-\frac{r_c}{\sin(\frac{\alpha}{2})}, 0.5 \right).\nonumber
\end{equation}
When $ \pi/2 \leq \alpha < \pi$, the left part of the boundary consists of three parts $\Gamma_4 = \Gamma_4^A \cup \Gamma_4^B \cup \Gamma_4^C$. The restriction that the arc fits into the unit square takes the form
\begin{equation}
0 < r_c < \frac{0.5}{\cos{(\frac{\alpha}{2})}}\nonumber
\end{equation}
 and the center of the arc is at point 
\begin{equation}
C = \left(\frac{0.5}{\tan(\frac{\alpha}{2})} - \frac{r_c}{\sin(\frac{\alpha}{2})}, 0.5\right).\nonumber
\end{equation}
As we use SI units the edge length of the square is $\SI{1}{\meter}$.

\subsection{Boundary conditions}
We impose the same boundary conditions on all computational domains. The domain is loaded by the shearing force $F=\SI{100}{\mega \pascal}$ that acts on the top and on the bottom edge of the unit square. Therefore $g=F$ on $\Gamma_3$ (top), $g=-F$ on $\Gamma_1$ (down) and $g=0$ on the remaining parts of the boundary. Boundary conditions \eqref{eq4:findingbc} take the form
\begin{align}
A_0^{g}(x_1, x_2) = \begin{cases}
-F x_1 & \mbox{on } \Gamma_1,\\
-F & \mbox{on } \Gamma_2,\\
-F x_1 & \mbox{on } \Gamma_3,\\
0 & \mbox{on } \Gamma_4,
\end{cases}
\nonumber
\end{align}
where $F=\SI{1e8}{\pascal}$.

\subsection{Variational problem on the space of finite elements}
Now we formulate the boundary value problem documented in \defref{d05:antiPlaneStressAiryFunction} as a variational problem on the space of finite elements. To derive the weak formulation, we multiply the first equation in \eqref{eq5:airybvp} by a test function $\phi \in W^{1,q'}_{0}(B)$, integrate over $B$ and use integration by parts. For the power-law solid \eqref{eq5:powerlaw}  where  $\sigma$ has the form
\begin{equation}
\sigma(|\T_v|) \sim  \left( 1 + |\T_v|^2\right)^{\frac{q'-2}{2}},\nonumber
\end{equation}
we seek the function $A \in W^{1,q'}(B)$ that satisfies the formulation
\begin{equation}
\begin{aligned}
\INTEGRAL{_B}{\sigma(\ABS{\grad A})\grad A \cdot \grad \phi}{x} &= 0  &&\forall \phi \in  W^{1,q'}_{0}(B),\\
A &= A_0^{g} &&\mbox{ on } \partial B,
\end{aligned}
\label{eq:weekformanti}
\end{equation}
where $A_0^g$  is specified in \eqref{eq4:findingbc}. When the response is given by Hooke's law \eqref{eq5:hookelaw}, where $\sigma$ is a constant, we set $q'=2$ and seek Airy's function in the space $A \in W^{1,2}(B)$.

To find numerical solutions, we first construct a triangulation $\mathcal{T}_h$ of the computational domain $B$. We employ a finite dimensional space $V_h\subset W^{1,q'}(B)$ over $\mathcal{T}_h$ and its counterpart with zero trace 
\begin{equation}
V_h^0 = \{ u \in V_h,  u|_{\partial B} = 0\}. \nonumber
\end{equation}
\begin{definition}[Airy's solution to the discrete Problem (P) within the context of the anti-plane stress setting]
\label{d05:antiPlaneStressAiryFem}
Let $B \subset \mathbb{R}^2$ be a computational domain described in \figref{fg5:computationalDomainVC}. Let $\mathcal{T}_h$ be a triangulation over $\overline{B}$. Let the functions ($\sigma$, $g$)  represent the problem data. We say that the function $A_h \in V_h \subset W^{1,q'}(B)$ is an Airy's function that solves the discrete Problem (P) within the context of the anti-plane stress setting, if
\begin{equation}
\begin{aligned}
\sum_{T \in \mathcal{T}_h}
\INTEGRAL{_T}{\sigma(\ABS{\grad A_h})\grad A_h \cdot \grad \phi_h}{x} &= 0  &&\forall \phi_h \in  V_h^0,\\
A_h &= A_0^{g} &&\mbox{ on } \partial B,
\end{aligned}
\label{eq5:discformanti}
\end{equation}
where $A_0^g$  is the function defined in \eqref{eq4:findingbc}.
\end{definition}

As $V_h$ we use the space of piece-wise continuous second order polynomials $V_h = X_h^2$, where
\begin{equation}
X_h^2  = \{ u_h \subset C^0(\overline{B}), u_h|_{T} \in P^2({T}), \forall T \in \mathcal{T}_h \}. \nonumber
\end{equation} 
We use Langrange elements of second order so that the $X_h^2 \subseteq W^{1,q'}(B)$ and therefore we can omit the sum over elements in \eqref{eq5:discformanti}. For the problems of the type \eqref{eq5:discformanti}, the quasi-norm interpolation error estimates were established, see for example \cite{ebmeyer_quasi_2005}.

For solving the problem computationally, we have developed software in Python that utilizes the FEniCS library, see \cite{alnaes_fenics_2015}. To find a solution for nonlinear power-law problems, we use damped Newton method, where the convergence criterion was based on the $L^2$ norm of the residuum. Computationally intensive tasks were performed using a cluster infrastructure supported by the Charles university, see \url{http://cluster.karlin.mff.cuni.cz/}.

\subsection{Material parameters used for simulations of titanium alloys}
\label{parameters}
We study the response of three different titanium alloys, namely Gum Metal, \textit{Ti-30Nb-10Ta-5Zr} alloy and \textit{Ti-24Nb-4Zr-7.9Sn} alloy. For each alloy, we compare four different models derived in \cite{kulvait_modeling_2017}. The parameters were estimated based on the general model \eqref{eq2:powerlaw1} and subsequently they are applied to the model \eqref{eq5:powerlaw}. The model labeled NLB is the power law model with parameters obtained as the best fit under the nonlinear bulk response condition and the model NLS is a model with parameters obtained under the nonlinear shear response condition, see \cite{kulvait_modeling_2017}\footnote{We are not using the same model and obtaining two different sets of values for the material parameters under two loading conditions. The NLB and NLS models are two different models that provide equally good fits. We used different experiments to obtain the material parameters for the two different models. These different sets of values fit the data equally well for the both the experiments. Which of these two models better explains the response of the model can only be determined by having other experiments against which these models can be corroborated. Such experiments are not available at the present time.}. The model LIN is a linear model \eqref{eq5:hookelaw}, where the shear modulus is given by Voigt-Reuss-Hill approximation. The original model \eqref{eq2:powerlaw1} has the parameters $(\tau_0, s', K, q', \mu)$ that appear in its representation. In the anti-plane stress setting, the parameters $s'$ and $K$ do not enter the model \eqref{eq5:powerlaw}, and therefore we work with a reduced set of parameters $(\tau_0, q', \mu)$. Parameters of computations are summarized in \tabref{tab05:comppar}.

\begin{table}[ht]
\centering
\begin{tabular}{l|lSSSS}
 Material & Model & {$\mu_L [\si{\giga \pascal}]$} & {$\tau_0 [\si{\giga \pascal}]$} & {$q'$} & {$\mu [\si{\giga \pascal}]$} \\ 
  \hline
\hline
Gum Metal & LIN1 & 23.5 & {-} & {-} & {-} \\ 
  \iftoggle{NLO}{Gum Metal & NLO1 & {-} & 0.5 & 1.92 & 20.9 \\}{}
  Gum Metal & NLB1 & {-} & 0.5 & 2.23 & 20.2 \\ 
  Gum Metal & NLS1 & {-} & 0.5 & 7.65 & 18668 \\ 
  \textit{Ti-30Nb-10Ta-5Zr} & LIN2 & 21.75 & {-} & {-} & {-} \\ 
  \iftoggle{NLO}{\textit{Ti-30Nb-10Ta-5Zr} & NLO2 & {-} & 0.5 & 1.88 & 24.5 \\}{}
  \textit{Ti-30Nb-10Ta-5Zr} & NLB2 & {-} & 0.5 & 2.49 & 22.3 \\ 
  \textit{Ti-30Nb-10Ta-5Zr} & NLS2 & {-} & 0.5 & 9.15 & 1001 \\ 
  \textit{Ti-24Nb-4Zr-7.9Sn} & LIN3 & 22.05 & {-} & {-} & {-} \\ 
  \iftoggle{NLO}{\textit{Ti-24Nb-4Zr-7.9Sn} & NLO3 & {-} & 0.5 & 2.14 & 18.6 \\}{}
  \textit{Ti-24Nb-4Zr-7.9Sn} & NLB3 & {-} & 0.5 & 2.99 & 16.5 \\ 
  \textit{Ti-24Nb-4Zr-7.9Sn} & NLS3 & {-} & 0.5 & 15.68 & 3378 \\ 
  \end{tabular}
\caption[Parameters of computer simulations]{Parameters of computer simulations. Linear model LIN with the response \eqref{eq5:hookelaw} is parametrized by the shear modulus $\mu_L$. Power law models \iftoggle{NLO}{NLB, NLS and NLO}{NLB and NLS} with the response \eqref{eq5:powerlaw} are parametrized by $\tau_0$, $q'$ and $\mu$. Data for the models are based on previous experimental corroboration, see \cite{kulvait_modeling_2017}.} 
\label{tab05:comppar}
\end{table}
  \section{Results}

\begin{figure}
\centering \includegraphics[width=\textwidth,height=\textheight,keepaspectratio]{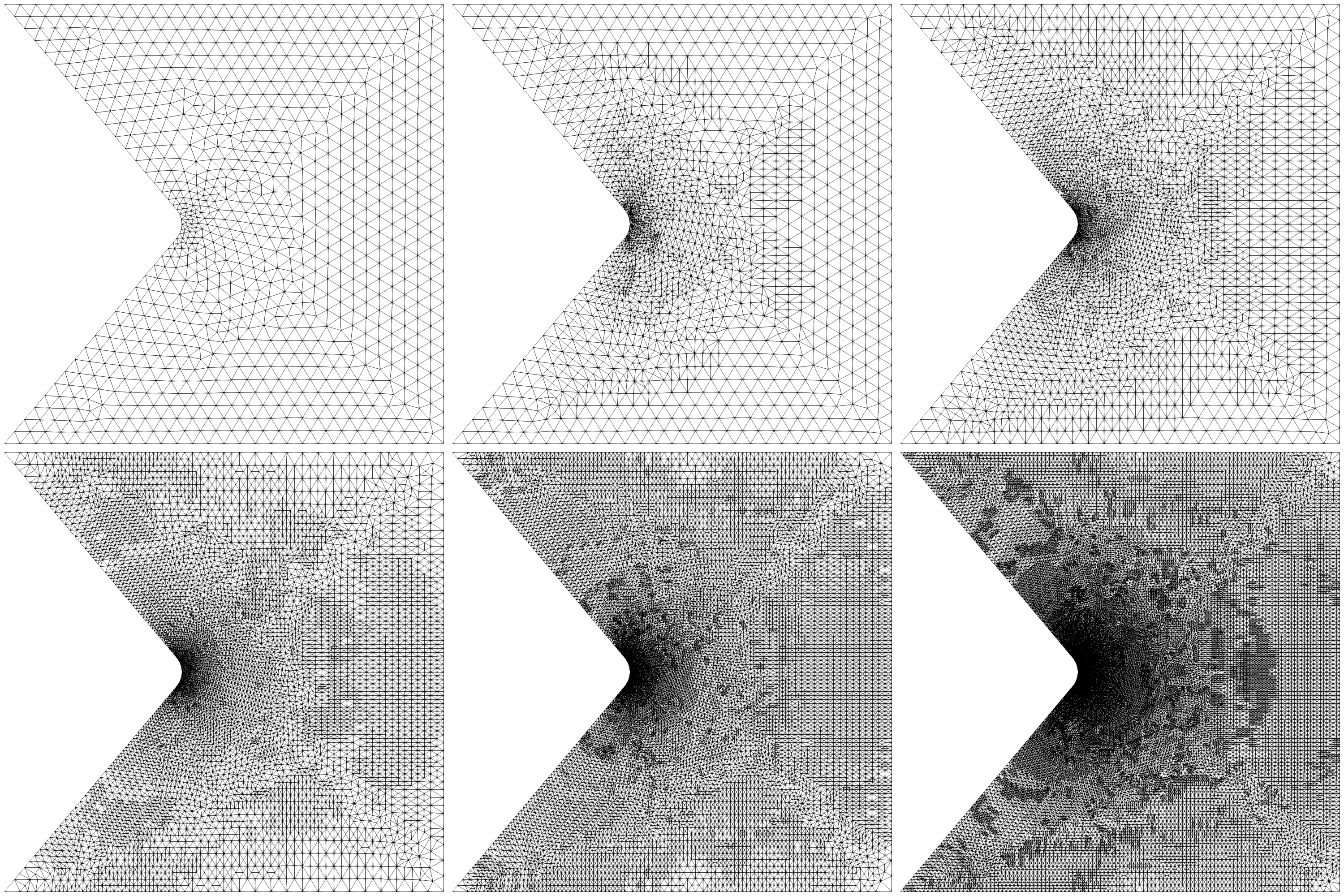}
    \caption[Visualization of adaptive mesh refinement, VC geometry]{Visualization of the adaptive refinement of the mesh in 5 steps. In upper left, there is the unrefined  mesh and in lower right, there is the mesh after 5 refinemet steps.  Visualisation for VC geometry, $\alpha = 100^\circ$, $r_c = 0.05$. }
    \label{img5:meshVC100}
\end{figure}

We have performed the following set of simulations. We computed solutions for $\alpha \in  \{\ang{1}, \ang{2}, \ldots \ang{70}, \ang{80} \ldots \ang{180} \}$. For each angle $\alpha$, we considered $r_c \in \{\num{0.05}, \num{0.01}, \num{0.005}, \num{0.001}\}$. In total, we use \num{324} different geometries of the computational domain. For each geometry, we generate one basic mesh and use $5$ adaptive refinements. In total, we use \num{1944} different meshes. Meshes were adaptively refined according to the local error estimator for the linear problem, this process is illustrated in Figure \ref{img5:meshVC100}. Mesh construction and refinement were performed using COMSOL Multiphysics software, version $3.5a$, see \cite{comsol_2008}. On each mesh, we performed computer simulations for \iftoggle{NLO}{\num{12}}{\num{9}} different model settings described in \tabref{tab05:comppar}, see \secref{parameters}.

 \subsection{Global convergence of solutions}
\label{sc5:gcs}
The primary goal of computer simulations is to study distributions of strains and stresses for each model and to find differences between the linear and nonlinear solutions and between solutions of the problems NLB and NLS that differ in the magnitude of the shear response exponent $q'$.

Prior to studying this, we investigate how accurate and reliable the solutions are. In this section, we study the global stability of the solutions with respect to a refinement level. To measure the relative error of solution with respect to the reference solution, we define the relative error norm.
\begin{definition}[Relative error norm]
Let $A_{ref}$ be a function (reference solution) on a finite dimensional space $V_h$. Let $A$ be a function (solution) on a finite dimensional space $Q_h \subseteq V_h$. Let $\NORM{.}$ denote a norm on $V_h$. Then the relative error norm of $A$ with respect to $A_{ref}$ is defined through
\begin{equation}
\NORM{A}^{rel}  = \frac{\NORM{A - A_{ref}}}{\NORM{A_{ref}}}.\label{eq5:relerrnrm}
\end{equation}
\end{definition}

For a sequence of adaptively refined meshes, the triangulation on one level is a subset of triangulations on the levels above. Let $V_h^i$ be a finite dimensional space of solutions on the $i$-th refinement level, then $V_h^0 \subseteq V_h^1 \ldots \subseteq V_h^5$. When analysing adaptively refined problems, we use a solution on the densest mesh, that is on the 5-th refinement level, as the reference solution. Then we compute the relative error norm for the solutions computed on croaser meshes. Our spaces $V_h^i$ consists of piece-wise continuous quadratic polynomials. Since $V_h^i \subseteq W^{1,2}(B)$ for $i \in \{0 \ldots 5\}$, we use the norm $\NORM{.} = \NORM{.}_{1,2}$ in \eqref{eq5:relerrnrm} for all problems.

\subsubsection*{Tables of mesh properties and error norms}
In Tables \ref{tb5:errornormsVC_alpha100_rc0.001}--\ref{tb5:errornormsVC_alpha100_rc0.05} we list some important parameters with regard to the computations with respect to the mesh refinement. In particular, we list the number of elements, the number of degrees of freedom, the mesh size parameters $h_{min}$ and $h_{max}$ and the relative error norms for each problem listed in \tabref{tab05:comppar}. For example, $\NORM{A_{NLB2}}^{rel}_{1,2}$ in a row in the table denotes the error norms with respect to the mesh refinement for the solution of the problem NLB2. In each table, the error norms for all the linear problems are compressed to the single row with a label $\NORM{A_{LIN}}^{rel}_{1,2}$, because for every mesh, the solution $A_{LIN}$ to the linear problem does not depend on the shear modulus $\mu_L$. It means that all the linear problems have identical Airy's functions and consequently error norms are identical as well. We note that
\begin{equation}
\NORM{A_{LIN}}^{rel}_{1,2} = \NORM{A_{LIN1}}^{rel}_{1,2} = \NORM{A_{LIN2}}^{rel}_{1,2} = \NORM{A_{LIN3}}^{rel}_{1,2}.
\end{equation}
We have chosen the following illustrative set of geometries to include in the tables. We consider data for $\alpha=\ang{100}$ and $r_c \in \{\num{0.001}, \num{0.01}, \num{0.05}\}$, see Tables \ref{tb5:errornormsVC_alpha100_rc0.001}--\ref{tb5:errornormsVC_alpha100_rc0.05}. 

\begin{table}[htbp]
\centering
\begin{tabular}{lllllll}
\toprule
 Refinement                    & 0              & 1              & 2              & 3               & 4               & 5              \\
\midrule
 Elements                      & $\num{2504}$   & $\num{6231}$   & $\num{13442}$  & $\num{28305}$   & $\num{58294}$   & $\num{127106}$ \\
 DOFs                          & $\num{5169}$   & $\num{12644}$  & $\num{27159}$  & $\num{57066}$   & $\num{117185}$  & $\num{255155}$ \\
 $h_{min}$                     & $\num{0.0002}$ & $\num{0.0002}$ & $\num{0.0002}$ & $\num{0.0002}$  & $\num{8e-05}$   & $\num{6e-05}$  \\
 $h_{max}$                     & $\num{0.04}$   & $\num{0.04}$   & $\num{0.04}$   & $\num{0.04}$    & $\num{0.03}$    & $\num{0.02}$   \\
 $\NORM{A_{LIN}}^{rel}_{1,2}$  & $\num{0.0006}$ & $\num{0.0003}$ & $\num{0.0002}$ & $\num{0.0002}$  & $\num{0.0002}$  & $\num{0.0}$    \\
 $\NORM{A_{NLB1}}^{rel}_{1,2}$ & $\num{0.0006}$ & $\num{0.0003}$ & $\num{0.0002}$ & $\num{0.0002}$  & $\num{0.0001}$  & $\num{0.0}$    \\
 $\NORM{A_{NLB2}}^{rel}_{1,2}$ & $\num{0.0005}$ & $\num{0.0003}$ & $\num{0.0002}$ & $\num{0.0002}$  & $\num{0.0001}$  & $\num{0.0}$    \\
 $\NORM{A_{NLB3}}^{rel}_{1,2}$ & $\num{0.0005}$ & $\num{0.0002}$ & $\num{0.0002}$ & $\num{0.0001}$  & $\num{0.0001}$  & $\num{0.0}$    \\
 \iftoggle{NLO}{
 $\NORM{A_{NLO1}}^{rel}_{1,2}$ & $\num{0.0006}$ & $\num{0.0003}$ & $\num{0.0002}$ & $\num{0.0002}$  & $\num{0.0002}$  & $\num{0.0}$    \\
 $\NORM{A_{NLO2}}^{rel}_{1,2}$ & $\num{0.0006}$ & $\num{0.0003}$ & $\num{0.0003}$ & $\num{0.0002}$  & $\num{0.0002}$  & $\num{0.0}$    \\
 $\NORM{A_{NLO3}}^{rel}_{1,2}$ & $\num{0.0006}$ & $\num{0.0003}$ & $\num{0.0002}$ & $\num{0.0002}$  & $\num{0.0001}$  & $\num{0.0}$    \\}{}
 $\NORM{A_{NLS1}}^{rel}_{1,2}$ & $\num{0.0004}$ & $\num{0.0002}$ & $\num{0.0001}$ & $\num{0.00008}$ & $\num{0.00006}$ & $\num{0.0}$    \\
 $\NORM{A_{NLS2}}^{rel}_{1,2}$ & $\num{0.0004}$ & $\num{0.0002}$ & $\num{0.0001}$ & $\num{0.00008}$ & $\num{0.00006}$ & $\num{0.0}$    \\
 $\NORM{A_{NLS3}}^{rel}_{1,2}$ & $\num{0.0004}$ & $\num{0.0002}$ & $\num{0.0001}$ & $\num{0.00008}$ & $\num{0.00005}$ & $\num{0.0}$    \\
\bottomrule
\end{tabular}\caption[Mesh properties and error norms, VC domain $\alpha = \ang{2}$, $r_c = \num{0.001}$]{Mesh properties and error norms of the solutions with respect to the refinement level for $\alpha = \ang{100}$ and $r_c = \num{0.001} $. For detailed description of row labels, see \secref{sc5:gcs} and \tabref{tab05:comppar}.}
\label{tb5:errornormsVC_alpha100_rc0.001}
\end{table}

\begin{table}[htbp]
\centering
\begin{tabular}{lllllll}
\toprule
 Refinement                    & 0              & 1              & 2              & 3              & 4               & 5              \\
\midrule
 Elements                      & $\num{2028}$   & $\num{4681}$   & $\num{9667}$   & $\num{20025}$  & $\num{40747}$   & $\num{87058}$  \\
 DOFs                          & $\num{4197}$   & $\num{9518}$   & $\num{19550}$  & $\num{40428}$  & $\num{82038}$   & $\num{174841}$ \\
 $h_{min}$                     & $\num{0.002}$  & $\num{0.001}$  & $\num{0.001}$  & $\num{0.0006}$ & $\num{0.0003}$  & $\num{0.0002}$ \\
 $h_{max}$                     & $\num{0.04}$   & $\num{0.04}$   & $\num{0.04}$   & $\num{0.04}$   & $\num{0.04}$    & $\num{0.02}$   \\
 $\NORM{A_{LIN}}^{rel}_{1,2}$  & $\num{0.001}$  & $\num{0.001}$  & $\num{0.0007}$ & $\num{0.0003}$ & $\num{0.0001}$  & $\num{0.0}$    \\
 $\NORM{A_{NLB1}}^{rel}_{1,2}$ & $\num{0.001}$  & $\num{0.0009}$ & $\num{0.0006}$ & $\num{0.0002}$ & $\num{0.00009}$ & $\num{0.0}$    \\
 $\NORM{A_{NLB2}}^{rel}_{1,2}$ & $\num{0.001}$  & $\num{0.0008}$ & $\num{0.0006}$ & $\num{0.0002}$ & $\num{0.00009}$ & $\num{0.0}$    \\
 $\NORM{A_{NLB3}}^{rel}_{1,2}$ & $\num{0.0009}$ & $\num{0.0007}$ & $\num{0.0005}$ & $\num{0.0002}$ & $\num{0.00008}$ & $\num{0.0}$    \\
 \iftoggle{NLO}{
 $\NORM{A_{NLO1}}^{rel}_{1,2}$ & $\num{0.001}$  & $\num{0.001}$  & $\num{0.0007}$ & $\num{0.0003}$ & $\num{0.0001}$  & $\num{0.0}$    \\
 $\NORM{A_{NLO2}}^{rel}_{1,2}$ & $\num{0.001}$  & $\num{0.001}$  & $\num{0.0007}$ & $\num{0.0003}$ & $\num{0.0001}$  & $\num{0.0}$    \\
 $\NORM{A_{NLO3}}^{rel}_{1,2}$ & $\num{0.001}$  & $\num{0.0009}$ & $\num{0.0006}$ & $\num{0.0002}$ & $\num{0.0001}$  & $\num{0.0}$    \\}{}
 $\NORM{A_{NLS1}}^{rel}_{1,2}$ & $\num{0.0007}$ & $\num{0.0005}$ & $\num{0.0004}$ & $\num{0.0001}$ & $\num{0.00005}$ & $\num{0.0}$    \\
 $\NORM{A_{NLS2}}^{rel}_{1,2}$ & $\num{0.0006}$ & $\num{0.0005}$ & $\num{0.0003}$ & $\num{0.0001}$ & $\num{0.00005}$ & $\num{0.0}$    \\
 $\NORM{A_{NLS3}}^{rel}_{1,2}$ & $\num{0.0006}$ & $\num{0.0005}$ & $\num{0.0003}$ & $\num{0.0001}$ & $\num{0.00005}$ & $\num{0.0}$    \\
\bottomrule
\end{tabular}\caption[Mesh properties and error norms, VC domain $\alpha = \ang{2}$, $r_c = \num{0.01}$]{Mesh properties and error norms of the solutions with respect to the refinement level for $\alpha = \ang{100}$ and $r_c = 0.01$. For detailed description of row labels, see \secref{sc5:gcs} and \tabref{tab05:comppar}.}
\label{tb5:errornormsVC_alpha100_rc0.01}
\end{table}

\begin{table}[htbp]
\centering
\begin{tabular}{lllllll}
\toprule
 Refinement                    & 0             & 1              & 2              & 3              & 4              & 5              \\
\midrule
 Elements                      & $\num{1758}$  & $\num{4153}$   & $\num{8481}$   & $\num{17287}$  & $\num{34819}$  & $\num{70086}$  \\
 DOFs                          & $\num{3645}$  & $\num{8460}$   & $\num{17172}$  & $\num{34920}$  & $\num{70138}$  & $\num{140827}$ \\
 $h_{min}$                     & $\num{0.01}$  & $\num{0.002}$  & $\num{0.002}$  & $\num{0.001}$  & $\num{0.0008}$ & $\num{0.0006}$ \\
 $h_{max}$                     & $\num{0.04}$  & $\num{0.04}$   & $\num{0.04}$   & $\num{0.04}$   & $\num{0.04}$   & $\num{0.03}$   \\
 $\NORM{A_{LIN}}^{rel}_{1,2}$  & $\num{0.003}$ & $\num{0.001}$  & $\num{0.0004}$ & $\num{0.0003}$ & $\num{0.0002}$ & $\num{0.0}$    \\
 $\NORM{A_{NLB1}}^{rel}_{1,2}$ & $\num{0.003}$ & $\num{0.0009}$ & $\num{0.0004}$ & $\num{0.0003}$ & $\num{0.0002}$ & $\num{0.0}$    \\
 $\NORM{A_{NLB2}}^{rel}_{1,2}$ & $\num{0.003}$ & $\num{0.0009}$ & $\num{0.0003}$ & $\num{0.0003}$ & $\num{0.0002}$ & $\num{0.0}$    \\
 $\NORM{A_{NLB3}}^{rel}_{1,2}$ & $\num{0.003}$ & $\num{0.0008}$ & $\num{0.0003}$ & $\num{0.0003}$ & $\num{0.0002}$ & $\num{0.0}$    \\
\iftoggle{NLO}{
 $\NORM{A_{NLO1}}^{rel}_{1,2}$ & $\num{0.003}$ & $\num{0.001}$  & $\num{0.0004}$ & $\num{0.0003}$ & $\num{0.0002}$ & $\num{0.0}$    \\
 $\NORM{A_{NLO2}}^{rel}_{1,2}$ & $\num{0.003}$ & $\num{0.001}$  & $\num{0.0004}$ & $\num{0.0003}$ & $\num{0.0002}$ & $\num{0.0}$    \\
 $\NORM{A_{NLO3}}^{rel}_{1,2}$ & $\num{0.003}$ & $\num{0.0009}$ & $\num{0.0004}$ & $\num{0.0003}$ & $\num{0.0002}$ & $\num{0.0}$    \\}{}
 $\NORM{A_{NLS1}}^{rel}_{1,2}$ & $\num{0.002}$ & $\num{0.0006}$ & $\num{0.0002}$ & $\num{0.0002}$ & $\num{0.0001}$ & $\num{0.0}$    \\
 $\NORM{A_{NLS2}}^{rel}_{1,2}$ & $\num{0.002}$ & $\num{0.0006}$ & $\num{0.0002}$ & $\num{0.0002}$ & $\num{0.0001}$ & $\num{0.0}$    \\
 $\NORM{A_{NLS3}}^{rel}_{1,2}$ & $\num{0.002}$ & $\num{0.0006}$ & $\num{0.0002}$ & $\num{0.0002}$ & $\num{0.0001}$ & $\num{0.0}$    \\
\bottomrule
\end{tabular}\caption[Mesh properties and error norms, VC domain $\alpha = \ang{2}$, $r_c = \num{0.05}$]{Mesh properties and error norms of the solutions with respect to the refinement level for $\alpha = \ang{100}$ and $r_c = \num{0.05}$. For detailed description of row labels, see \secref{sc5:gcs} and \tabref{tab05:comppar}.}
\label{tb5:errornormsVC_alpha100_rc0.05}
\end{table}
  \subsection{Comparison of solutions}
Now, we focus on the geometry that has been smoothened by a circle or an arc of diameter $r_c$. We consider $\alpha=\ang{90}$ and compare the behavior of the models LIN2, \iftoggle{NLO}{NLO2,} NLB2 and NLS2.

\begin{figure}
\vspace*{-1cm}
\centering \includegraphics[width=\textwidth,height=\textheight,keepaspectratio]{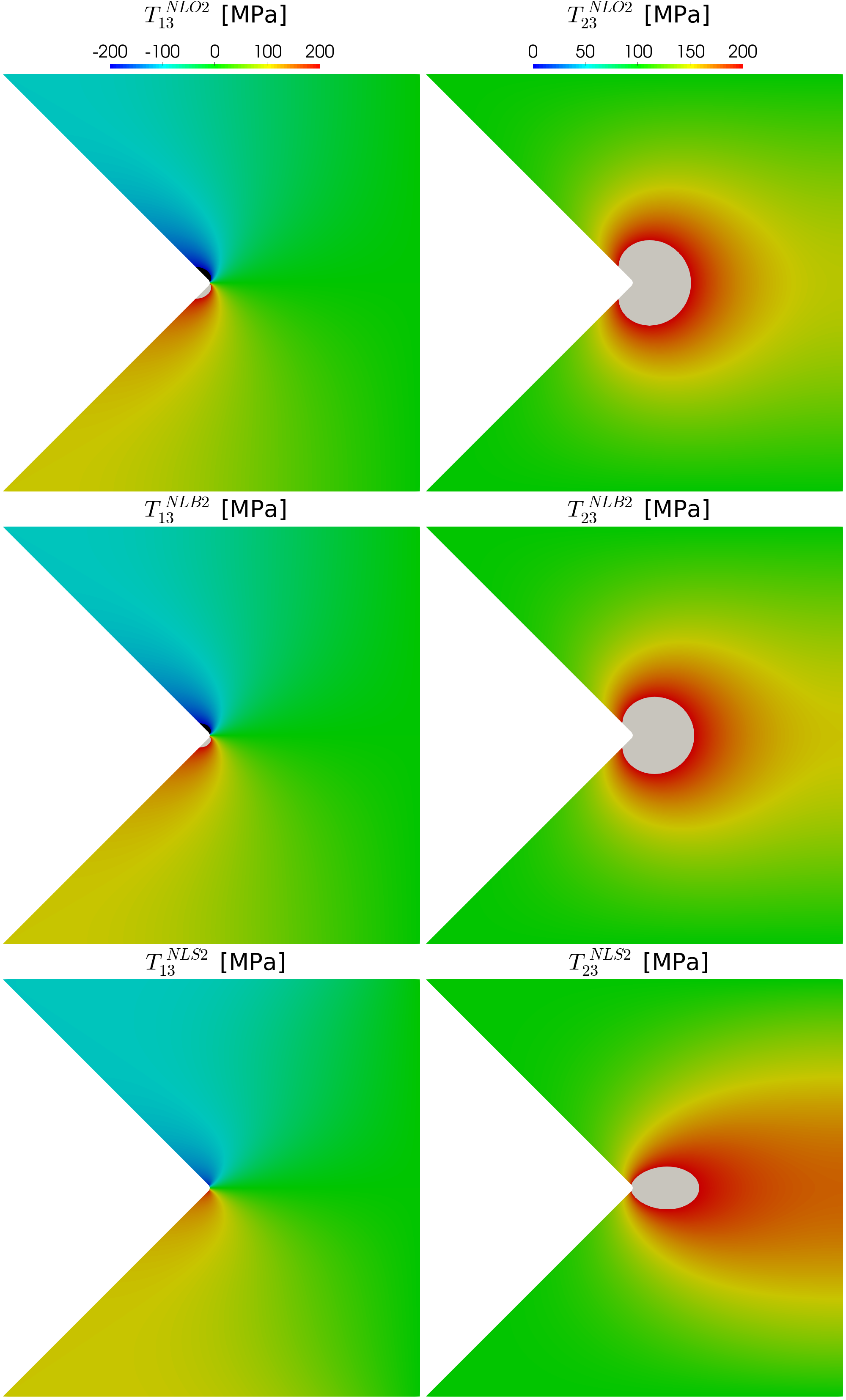}
        \caption[Comparison of stress distributions in VC geometry]{Comparison of stress distributions in VC geometry with $\alpha = \ang{90}$ and $r_c = \num{0.01}$ for models NLO2, NLB2, NLS2.}
    \label{img5:comparisonTVC90}
\end{figure}

\begin{figure}
\vspace*{-1cm}
\centering \includegraphics[width=\textwidth,height=\textheight,keepaspectratio]{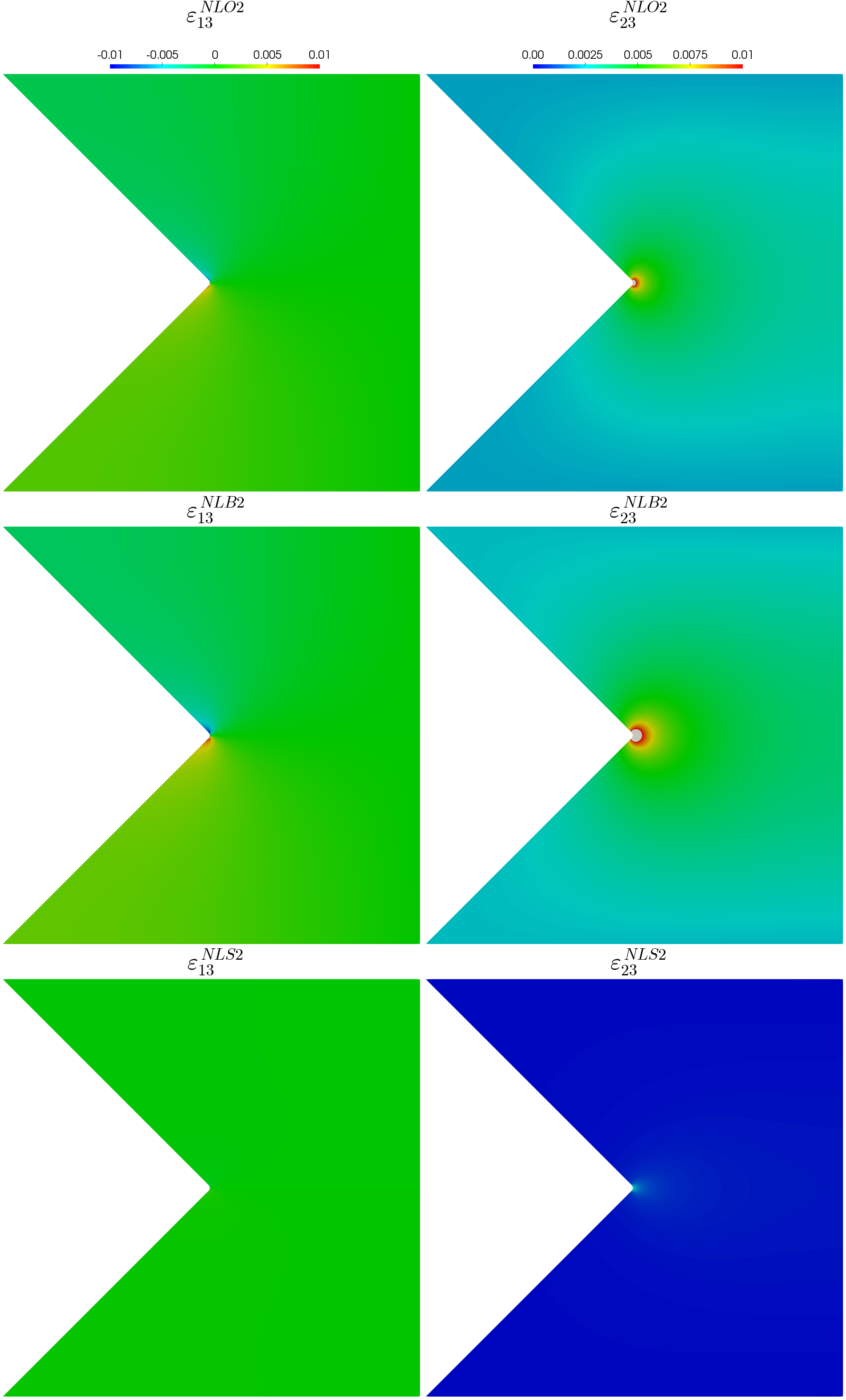}
        \caption[Comparison of strain distributions in VC geometry]{Comparison of strain distributions in VC geometry with $\alpha = \ang{90}$ and $r_c = \num{0.01}$ for models NLO2, NLB2, NLS2.}
    \label{img5:comparisonEVC90}
\end{figure}

\begin{figure}
\centering \includegraphics[width=\textwidth,height=\textheight,keepaspectratio]{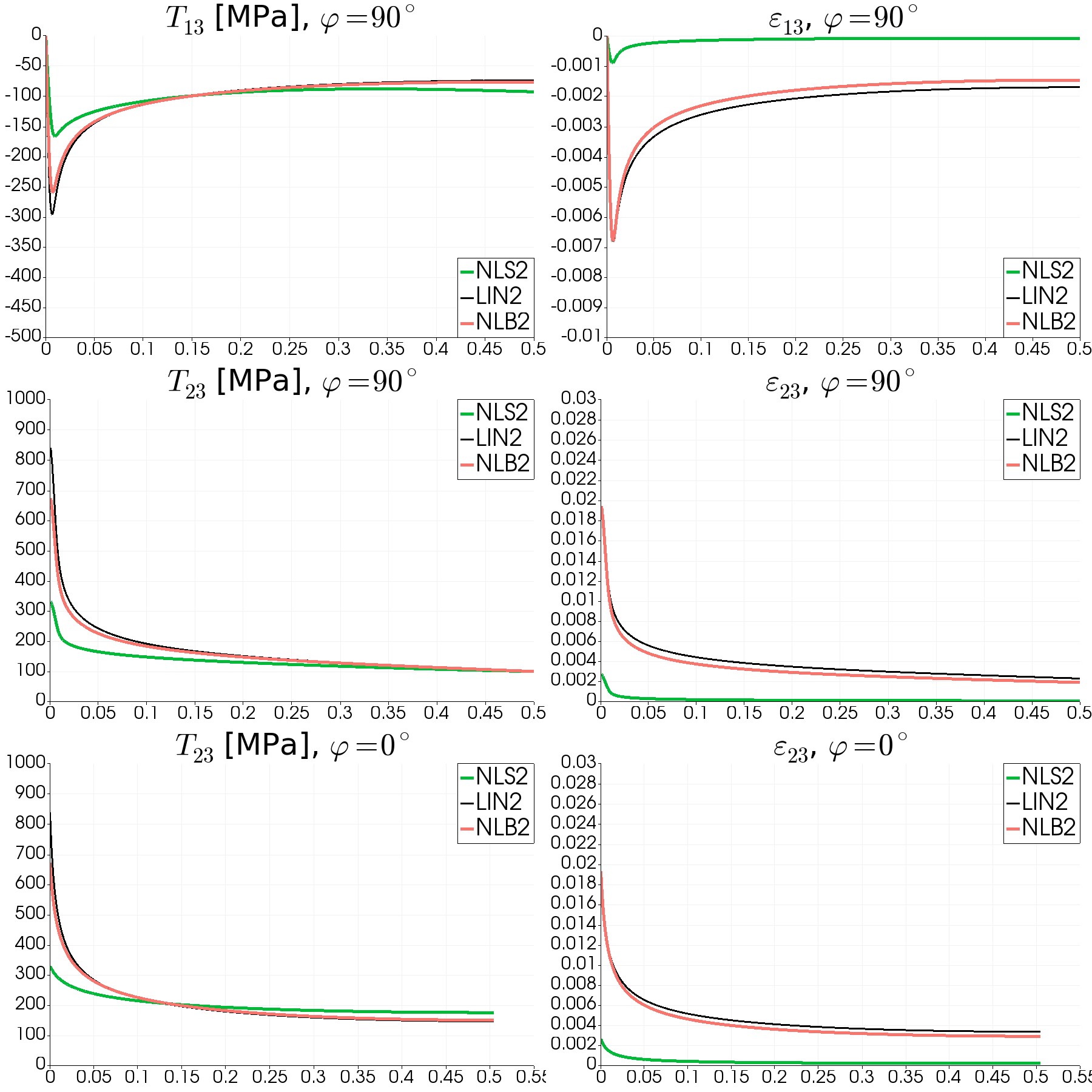}
    \caption[Plots of stress and strain over the line in VC geometry]{Plots of stress and strain for models LIN2, NLO2, NLB2 and NLS2 in VC geometry with $\alpha = \ang{90}$ and $r_c=\num{0.01}$. The plots are created from the line starting at the point $C + (r_c,0)$ of the geometry and ending either at point $(1,0.5)$, which is denoted $\varphi = \ang{0}$, or at point $C + (r_c,0.5)$, which corresponds to $\varphi = \ang{90}$.}
    \label{img5:comparisonVC90plots}
\end{figure}

In \figref{fg5:computationalDomainVC}, we consider $r_c = \num{0.01}$. Visualizations of the stress and strain distributions are depicted in  \figref{img5:comparisonTVC90} and \figref{img5:comparisonEVC90} respectively. We construct plots of stress and strain over the line from the smothened tip, that is the point $C + (r_c,0)$ to the right, which is denoted as $\varphi = \ang{0}$ and to the top, which is denoted as $\varphi = \ang{90}$, see \figref{img5:comparisonVC90plots}.

We find that there is no singularity in the stresses or the strains for the various situations considered in the study.

 \subsection{Dependence on $\alpha$ and $r_c$}
Thus far, we have been comparing the solutions when the geometrical setting is fixed. In this section, we study how the solutions depend on the parameters of the geometry. We consider the behavior of solutions when the diameter $r_c$ is decreasing and when the opening angle $\alpha$ is increasing.

\begin{figure}
\centering \includegraphics[width=\textwidth,height=\textheight,keepaspectratio]{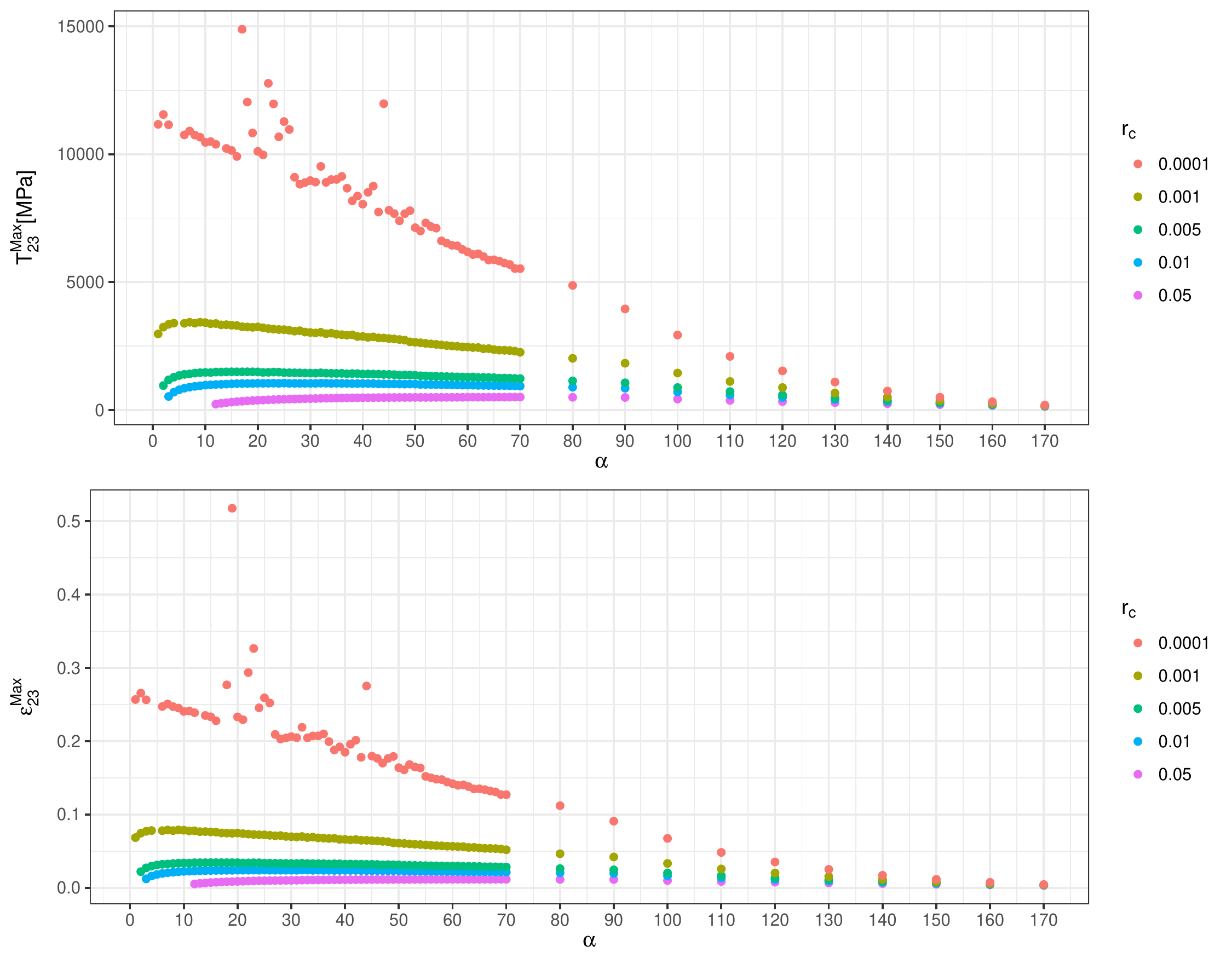}
    \caption[Maximal $\Tx{23}$ and $\Ex{23}$ in VC geometry, model LIN2]{Maximal $\Tx{23}$ and $\Ex{23}$ with respect to $\alpha$ for model LIN2.}
    \label{img5:VCLMAX}
\end{figure}

\begin{figure}
\centering \includegraphics[width=\textwidth,height=\textheight,keepaspectratio]{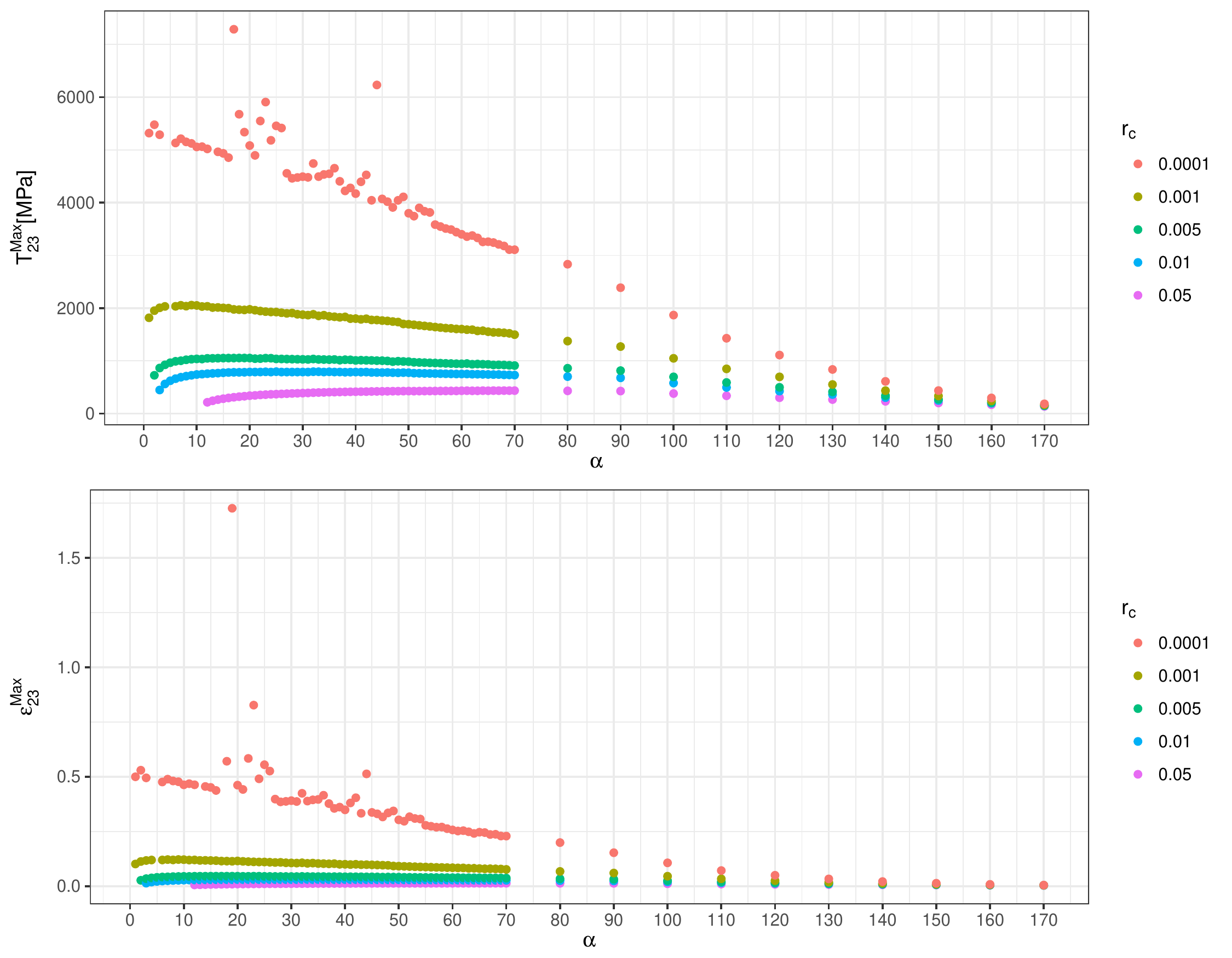}
    \caption[Maximal $\Tx{23}$ and $\Ex{23}$ in VC geometry, model NLB2]{Maximal $\Tx{23}$ and $\Ex{23}$ with respect to $\alpha$ for model NLB2.}
    \label{img5:VCBMAX}
\end{figure}

\begin{figure}
\centering \includegraphics[width=\textwidth,height=\textheight,keepaspectratio]{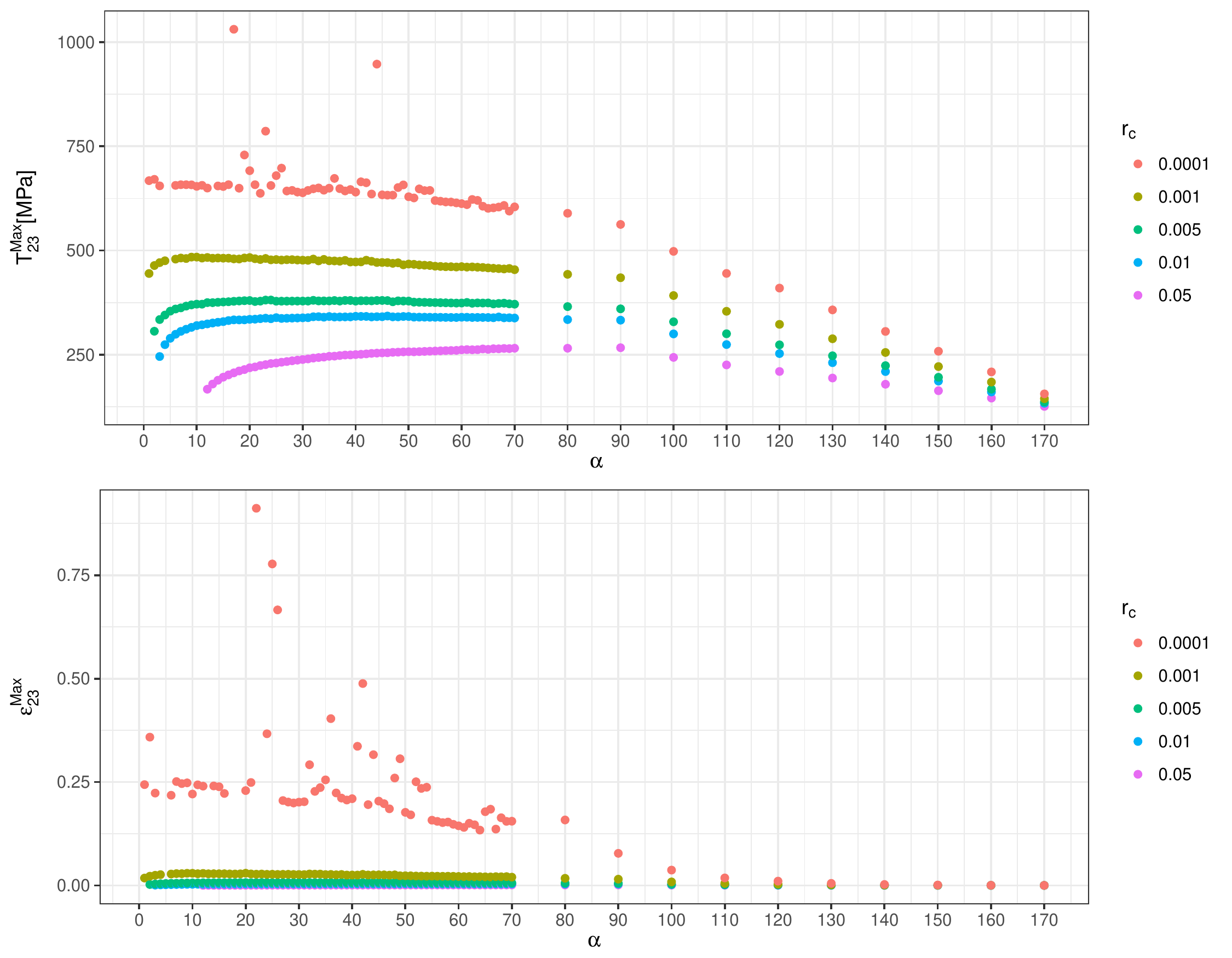}
    \caption[Maximal $\Tx{23}$ and $\Ex{23}$ in VC geometry, model NLS2]{Maximal $\Tx{23}$ and $\Ex{23}$ with respect to $\alpha$ for model NLS2.}
    \label{img5:VCSMAX}
\end{figure}

In Figures \ref{img5:VCLMAX}--\ref{img5:VCSMAX}, we plot the dependence of the highest value of $\Tx{23}$ and $\Ex{23}$ in the whole computational domain with respect to $\alpha$ for all the models that were studied. Plots for LIN2 model in \figref{img5:VCLMAX}. Plots for NLB2 and NLS2 are provided, respectively in Figures \ref{img5:VCBMAX}--\ref{img5:VCSMAX}.

Next we define the distance from the notch tip to the point, where $\Tx{23}$ or $\Ex{23}$ acquires a particular value.

\begin{definition}[Functions $dst^T$ and $dst^e$]
Let $(r, \varphi)$ be the cylindrical coordinate system centered at point $C+(r_c,0)$.  We define the functions $dst^T : [0, 2\pi) \times \mathbb{R} \to \mathbb{R}^+_0$ and $dst^e : [0, 2\pi) \times \mathbb{R} \to \mathbb{R}^+_0$ as follows
\begin{equation}
\begin{aligned}
dst^T(\varphi, T) &= 
\begin{cases} 
\infty &\mbox{if } \forall r \geq 0 : T < \Tx{23}(r, \varphi),\\
-\infty &\mbox{if } \forall r \geq 0 : T > \Tx{23}(r, \varphi),\\
\min\{r \geq 0,  \Tx{23}(r, \varphi)=T\} &\mbox{otherwise}. 
\end{cases}\\
dst^e(\varphi, e) &=
\begin{cases} 
\infty &\mbox{if } \forall r \geq 0 : e < \Ex{23}(r, \varphi),\\
-\infty &\mbox{if } \forall r \geq 0 : e > \Ex{23}(r, \varphi),\\
\min\{r \geq 0, \Ex{23}(r, \varphi)=e\} & \mbox{otherwise}. 
\end{cases}
\end{aligned}\nonumber
\end{equation}
\end{definition}

For visualizations, we define the three following distances called DT2, DT3 and DE1 such that
\begin{equation}
\begin{aligned}
DT2 &= dst^T(0, \SI{200}{\mega \pascal}),\\
DT3 &= dst^T(0, \SI{300}{\mega \pascal}),\\
DE1 &= dst^e(0, 0.01).
\end{aligned}
\label{eq5:distances}
\end{equation}

\begin{figure}
\centering \includegraphics[width=\textwidth,height=\textheight,keepaspectratio]{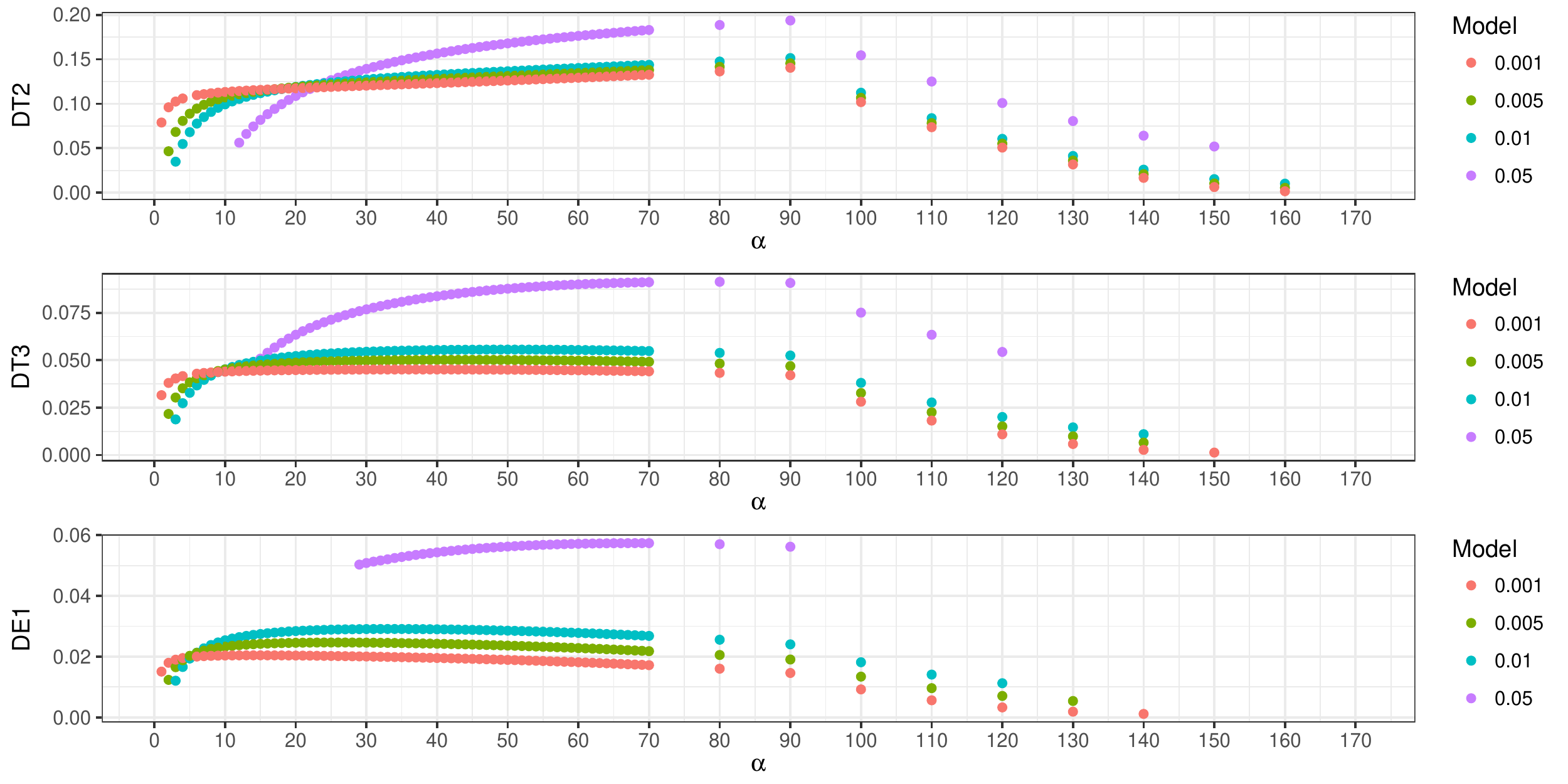}
    \caption[Distances DT2, DT3 and DE1, model LIN2]{Distances DT2, DT3 and DE1, see \eqref{eq5:distances}, with respect to $\alpha$ for model LIN2.}
    \label{img5:DCL}
\end{figure}

\begin{figure}
\centering \includegraphics[width=\textwidth,height=\textheight,keepaspectratio]{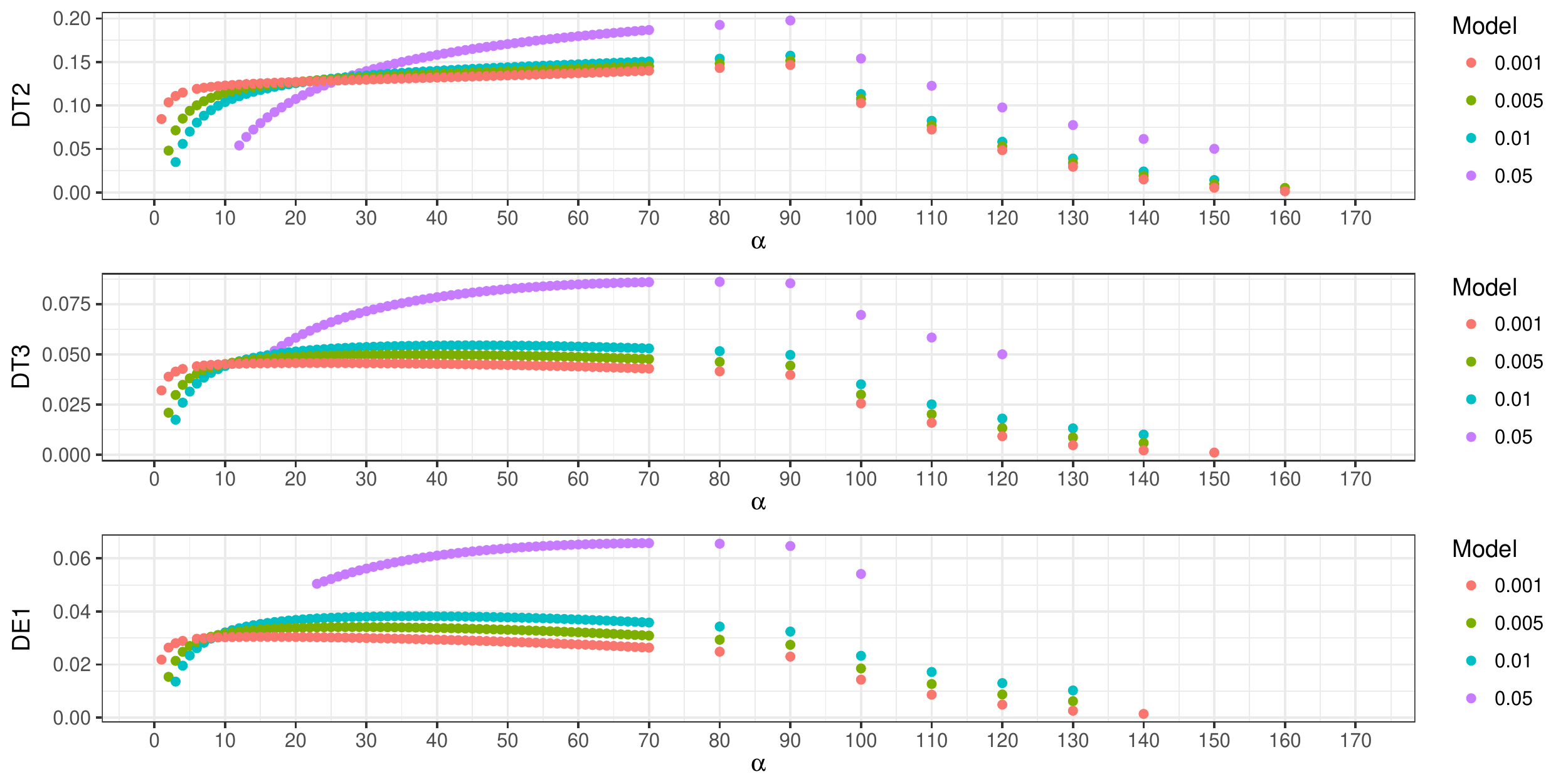}
    \caption[Distances DT2, DT3 and DE1 in VO and VC geometries, model NLB2]{Distances DT2, DT3 and DE1, see \eqref{eq5:distances}, with respect to $\alpha$ for model NLB2.}
    \label{img5:DCB}
\end{figure}

\begin{figure}
\centering \includegraphics[width=\textwidth,height=\textheight,keepaspectratio]{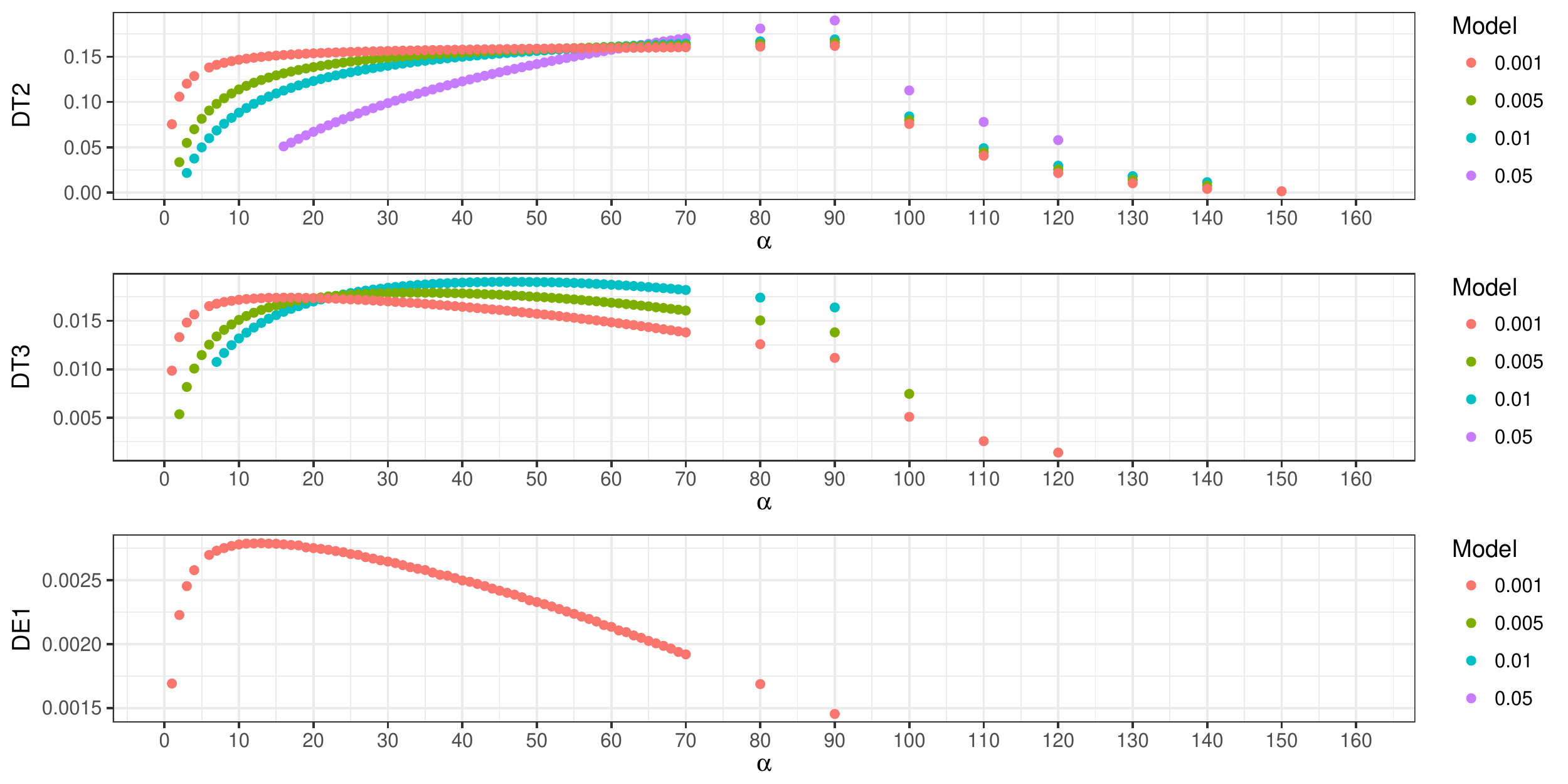}
    \caption[Distances DT2, DT3 and DE1 in VO and VC geometries, model NLS2]{Distances DT2, DT3 and DE1, see \eqref{eq5:distances}, with respect to $\alpha$ for model NLS2.}
    \label{img5:DCS}
\end{figure}

These distances are visualized for model LIN2 in \figref{img5:DCL}, for model NLB2 in \figref{img5:DCB} and for model NLS2 in \figref{img5:DCS}.

 \section{Discussion}
We considered finite element approximations of the Problem (P) in the reduced geometrical setting of anti-plane stress. We simulated the response of three titanium alloys, described by the power-law models \eqref{eq5:powerlaw}, with the material moduli previously identified in \cite{kulvait_modeling_2017}, see \tabref{tab05:comppar}. Numerical simulations were performed on adaptively refined triangular meshes. 

We have verified that the computed variational solutions are stable. In Tables \ref{tb5:errornormsVC_alpha100_rc0.001}--\ref{tb5:errornormsVC_alpha100_rc0.05}, we list relative error norms \eqref{eq5:relerrnrm} with respect to the refinement level for all the problems studied for numerous different parameters of the geometry. The differences between the computed solutions on the two most refined meshes are typically of order lower than $\num{e-4}$. The global stability measured by the error norms for the power-law models is comparable to or even better than the global stability for the linear problem. There is no singularity of stress and strain observed, see Figure \ref{img5:comparisonTVC90} and \figref{img5:comparisonVC90plots}.

It is evident that with increasing angle $\alpha$, the solution is becoming even more stable due to the decrease of the severity of the stress concentration. We observe that the solutions are more stable for higher values of the power law exponents $q'$. Namely the error norms for the model NLS2 with the highest $q'$ are consistently smaller than the error norms for the NLB2 model with the lower value of $q'$.

We measured the local error of computed solutions with respect to the refinement by the norm of the difference in the stress. 

When studying the behavior of the solutions with respect to parameters of the geometry, see Figures \ref{img5:VCLMAX}--\ref{img5:DCS}, we found that there is substantial stability of the solution in comparisons to the solutions across completely different meshes. We can clearly identify patterns in Figures \ref{img5:VCLMAX}--\ref{img5:DCS} that arise with respect to the notch opening angle $\alpha$.

 \bibliographystyle{spmpsci}   

\end{document}